\documentclass[3p]{article}

\usepackage[utf8]{inputenc} 
\usepackage[T1]{fontenc}    
\usepackage{hyperref}       
\usepackage{url}            
\usepackage{booktabs}       
\usepackage{amsfonts}       
\usepackage{nicefrac}       
\usepackage{microtype}      
\usepackage{lipsum}
\usepackage{graphicx}
\usepackage{amsmath}        
\usepackage{algorithm}

\usepackage{algpseudocode}
\usepackage{amsthm}         

\usepackage{datetime}
\usepackage{hyperref}
\usepackage{subfig}

\usepackage{caption}

\usepackage{placeins}       

\usepackage[verbose=true,letterpaper]{geometry}
\AtBeginDocument{
  \newgeometry{
    textheight=9in,
    textwidth=6.5in,
    top=1in,
    headheight=14pt,
    headsep=25pt,
    footskip=30pt
  }
}

\hypersetup{
    colorlinks=true,
    citecolor=blue,
    linkcolor=blue,
    filecolor=magenta,      
    urlcolor=cyan,
    pdfborderstyle={/S/U/W 1} 

}

\graphicspath{{media/}}     

\newtheorem{remark}{Remark}
\numberwithin{remark}{section}

\def\keywordname{{\bfseries \emph Keywords}}%
\def\keywords#1{\par\addvspace\medskipamount{\rightskip=0pt plus1cm
\def\and{\ifhmode\unskip\nobreak\fi\ $\cdot$
}\noindent\keywordname\enspace\ignorespaces#1\par}}

\title{MGCNN: a learnable multigrid solver for sparse linear systems from PDEs on structured grids
}

\author{
    Yan Xie, Minrui Lv, Chensong Zhang \\
    LSEC, Academy of Mathematics and Systems Science, \\
    Chinese Academy of Sciences, Beijing, 100190, China; \\
    School of Mathematical Sciences, \\
    University of Chinese Academy of Sciences, Beijing, 100049, China \\
    \texttt{\{xieyan2021, lvminrui, zhangcs\}@lsec.cc.ac.cn} \\
}

\begin{document}
\maketitle

\begin{abstract}
    This paper presents a learnable solver tailored to iteratively solve sparse linear systems from discretized partial differential equations (PDEs). Unlike traditional approaches relying on specialized expertise, our solver streamlines the algorithm design process for a class of PDEs through training, which requires only training data of coefficient distributions.
    The proposed method is anchored by three core principles: (1) a multilevel hierarchy to promote rapid convergence, (2) adherence to linearity concerning the right-hand-side of equations, and (3) weights sharing across different levels to facilitate adaptability to various problem sizes. Built on these foundational principles and considering the similar computation pattern of the convolutional neural network (CNN) as multigrid components, we introduce a network adept at solving linear systems from PDEs with heterogeneous coefficients, discretized on structured grids. Notably, our proposed solver possesses the ability to generalize over right-hand-side terms, PDE coefficients, and grid sizes, thereby ensuring its training is purely offline.
    
    To evaluate its effectiveness, we train the solver on convection-diffusion equations featuring heterogeneous diffusion coefficients. The solver exhibits swift convergence to high accuracy over a range of grid sizes, extending from $31 \times 31$ to $4095 \times 4095$. Remarkably, our method outperforms the classical Geometric Multigrid (GMG) solver, demonstrating a speedup of approximately 3 to 8 times. Furthermore, our numerical investigation into the solver's capacity to generalize to untrained coefficient distributions reveals promising outcomes. Specifically, training the solver on a mixture of several coefficient distributions appears to significantly enhance its robustness when faced with previously unseen distributions. This suggests a potential avenue for improving the solver's robustness across diverse scenarios.
\end{abstract}

\keywords{learnable multigrid solver \and sparse linear systems \and  convolutional neural network (CNN) \and heterogeneous coefficient \and convection-diffusion equations}

\section{Introduction}
\label{sec:intro}

\subsection{Background}
Iterative methods play a pivotal role in solving linear systems arising from the discretization of partial differential equations (PDEs). Renowned for their low complexity, parallel scalability, and the flexibility to control accuracy during iterations, these methods are commonly preferred for large-scale problems. However, no universally efficient iterative method exists for all types of linear systems; their efficacy is intrinsically linked to the specific problem. The creation of an effective iterative method usually involves expert knowledge and significant trial-and-error. In contrast, the emergence of deep learning has revolutionized multiple sectors, notably computer vision~\cite{he2016deep, vaswani2017attention, goodfellow2014generative}, natural language processing~\cite{vaswani2017attention, mikolov2013distributed, brown2020language, kenton2019bert}, and bioinformatics~\cite{alipanahi2015predicting, jumper2021highly, gligorijevic2021structure}. This trend towards data-driven methodologies is reshaping problem-solving strategies, a shift that is increasingly influencing numerical methods for solving PDE problems, as evidenced by recent developments in the field.

Broadly speaking, deep learning has been applied to solve PDEs in two ways. The first approach involves integrating neural networks into traditional algorithms, replacing certain modules. These modules are often tailored to specific problems or based on empirical knowledge. The goal is to enhance the performance of classical algorithms. For instance, neural networks have been used to pinpoint discontinuity locations~\cite{ray2019detecting}, introduce artificial viscosity for stabilization~\cite{discacciati2020controlling}, improve multigrid methods~\cite{katrutsa2017deep, greenfeld2019learning}, and identify optimal initial values for nonlinear problems~\cite{huang2020int}.

The second approach to applying deep learning in PDE solving involves using neural networks in roles akin to traditional PDE solvers. A prominent stream within this category is the Physics-Informed Neural Networks (PINNs)~\cite{raissi2019physics}. PINNs essentially employ neural networks to approximate the solution of a PDE. This is achieved through minimizing a residual loss that is informed by physical laws, thus adhering to the constraints of various physical equations and boundary conditions. The strength of PINNs lies in their ability to tackle high-dimensional problems and manage complex geometries, owing to their mesh-free nature. This approach has been successfully applied in a range of scenarios, such as solving both forward and inverse problems~\cite{bar2021strong, moseley2020solving, song2022versatile, stanziola2021helmholtz, yu2018deep, lu2021deepxde, raissi2019physics}. Their solve phase lies in training parameters for the network to approximate the solution function.

In our study, we align with the second approach, focusing on a pure data-driven algorithm that meets the problem-specific nature of iterative methods. Rather than solving a problem during training like PINNs, we dive into another stream to develop a network capable of learning an operator to solve the problem during inference. While there has been research on operator learning for continuous PDEs~\cite{kissas2022learning, lu2021learning, patel2021physics}, our interest lies in creating a solver for discretized linear systems. There already exist solvers capable of solving the given problem in a single~\cite{wiecha2019deep, xue2020amortized, khoo2019switchnet, grementieri2022towards} or multiple~\cite{hsieh2018learning, pfaff2020learning, rizzuti2019learned, CHEN2022110996, azulay2022multigrid, lerer2023multigrid} iterative applications of a feed-forward network. 

\subsection{Related Works}
\label{sec:relatedwork}

Classical numerical solvers for PDEs typically utilize a multilevel structure to ensure quick convergence and scalability, especially for large problem sizes. Inspired by this, we aim to replicate a similar multilevel structure in our solver network. Within the current learning framework, the convolutional neural network (CNN)~\cite{lecun1998gradient} stands out as the most feasible network architecture for this purpose. Consequently, we have chosen the CNN as the main building block of our proposed solver. This decision also guides our focus to problems that are discretized on structured grids, aligning with the inherent capabilities of CNNs. 

This concept of introducing a multilevel structure into building neural networks has parallels in the field of image processing, particularly in architectures like U-net~\cite{ronneberger2015u}, which have effectively implemented similar strategies. In fact, the work by He et al.~\cite{he2019mgnet} elucidated the connection between convolutional neural networks (CNNs) and multigrid algorithms, paving the way for the creation of MgNet. In the realm of numerical PDE solvers, there have been notable efforts to develop CNN-based multigrid neural networks. For example, Chen et al.~\cite{CHEN2022110996} introduced Meta-MgNet, a model capable of adapting to varying problem parameters. Azulay et al.~\cite{azulay2022multigrid} and Lerer et al.~\cite{lerer2023multigrid} integrated neural networks with the classical shift Laplacian method; they trained these networks on multi-level smooth right-hand-side datasets, enabling them to act as preconditioners for the Flexible Generalized Minimal Residual (FGMRES) method. This approach was particularly effective in solving the Helmholtz equation on heterogeneous media.

Moreover, the success of Google's GraphCast~\cite{lam2023learning} initially encouraged us to explore Artificial Intelligence (AI) as a potential multigrid solver in the field of solving sparse linear systems from PDEs.

\subsection{Our Contributions}
\label{sec:contribution}

In our endeavor to create a learnable solver for sparse linear systems from discretized PDEs, we anchor our design on three key principles: (1) employing a multi-level hierarchical computation structure, (2) implementing the network linearly on the right-hand-side of equations, and (3) sharing weights across different levels. These principles enable our network to be effective, free of expertise, and applicable to various problem sizes.
For problems set on structured grids, we observe that Meta-MgNet~\cite{CHEN2022110996} aligns with the first two principles. However, it's not tailored for heterogeneous problems. 
Other studies, such as those by Azulay et al.~\cite{azulay2022multigrid} and Lerer et al.~\cite{lerer2023multigrid}, have developed neural networks for heterogeneous Helmholtz problems. These only satisfy the first principle and necessitate integrating the shift Laplacian method~\cite{oosterlee2009shifted} during both training and solving phases. Furthermore, to the best of our knowledge, there appears to be a gap in existing research concerning the third principle. This principle, the sharing of weights across levels, is crucial for our goal of scaling the network to tackle large problem sizes efficiently during training.

Adhering to our established principles, we develop a specific network trained to address linear systems from PDEs with heterogeneous coefficients, discretized on structured grids. Our learnable solver can generalize over right-hand-side terms, PDE coefficients, and grid sizes. We train our solver for the convection-diffusion equation with heterogeneous diffusion coefficients, and it showcases rapid convergence to high accuracy across a range of grid sizes, up to $4095 \times 4095$. Our method achieves a speedup of approximately 3 to 8 times compared to a Geometric Multigrid (GMG) solver on a GPU. This accomplishment underscores the potential of training neural networks in effectively resolving large, sparse linear systems derived from discretized PDEs, without the need for algorithm design and additional expert input.

\subsection{Paper Organization}
The rest of this paper is organized as follows: Section~\ref{sec:preliminaries} lays out the necessary preliminaries and fundamental concepts that form the basis of our study. The terminologies established in this section are consistently employed throughout the paper. Section~\ref{sec:methodology} offers an in-depth discussion of our proposed methodology. In Section~\ref{sec:experiments}, we conduct a series of experiments to demonstrate the effectiveness and reliability of our approach. The paper concludes with Section~\ref{sec:conclusion}, where we summarize our findings and contemplate future avenues of research. For additional insights and detailed information about training and coefficient distributions, readers are directed to Appendix~\ref{apd:training} and Appendix~\ref{apd:coeff_dist}. Notably, as the inaugural endeavor employing the three principles to construct a learnable multigrid solver, we provide extensive discussion in titled remarks, particularly within Section~\ref{sec:methodology}. We encourage readers to explore interested remarks to gain deeper understanding on specific aspects of our methodology.

\section{Preliminaries and Notations}
\label{sec:preliminaries}
In this section, we lay the groundwork for our study by introducing the core concepts and basic methodologies. We start with a discussion on the solver as a decoder, followed by an overview of iterative methods. We then delve into the multigrid algorithm, a multi-level iterative method for solving linear systems derived from certain PDE systems. Next, we explore the properties of convolutional neural networks and their variant, transposed convolutional neural networks. Additionally, we present the extension of the Residual Network~(ResNet) to a multigrid framework, an essential aspect of our methodology. These discussions are pivotal for appreciating the neural network architecture we utilize. We conclude this section with the discretization method for the convection-diffusion equation, which serves as our test problem. All these preliminaries are \textit{crucial} for understanding the methodology, experiments, and results presented in the subsequent sections.

\subsection{Solver as Decoder}
\label{sec:solverdecoder}
\begin{figure}[htbp]
\centering
\includegraphics[width=0.91\textwidth]{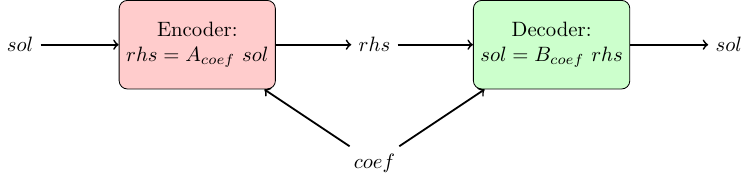}
\caption{Solver as Decoder}
\label{fig:solverdecoder}
\end{figure}

We first establish a broader perspective to clarify the problem and our task. Consider the linear system $A_{coef}\ sol = rhs$, where the linear operator or matrix defined by the coefficient tensor $coef$ acts as an encoder. This encoder transforms the solution $sol$ into the $rhs$. Our task is to design a decoder to reverse this process, decoding the information from $rhs$ to the solution $sol$. This decoder can be represented as a linear operator $B_{coef}$, i.e., $sol = B_{coef}\ rhs$. The process is illustrated in Fig.~\ref{fig:solverdecoder}. From this figure, it is important to note that our implementation of decoder $B_{coef}$ is independent of the specific format of encoder~$A_{coef}$. 

\subsection{Iterative Methods}
\label{sec:itermethods}
In the context of a linear system $A\ sol = rhs$, where $A$ is a square matrix, $sol$ is the solution vector, and $rhs$ is the right-hand-side vector, the objective of a linear solver is to find the solution $sol$ that meets satisfactory tolerance. We can obtain it by directly solving the equation or iteratively enhancing the solution's accuracy. The latter approach is often preferred due to its lower complexity and potential for typically high level of parallelization.

One such iterative method is the \textit{stationary iterative method}, which solves linear systems by iteratively addressing the error equation $A\ sol = r$, where $r = rhs - A\ sol$ represents the residual vector. This method can be defined as follows:
\begin{equation}
    \label{eq:stationary_iteration}
    sol^{k+1} = sol^k + B\ r^k = sol^k + B\ (rhs-A\ sol^k),
\end{equation}
where $sol^k$ is the solution at the $k$-th iteration and $B$ serves as the solver, providing an approximate error correction for each residual $r^k$. To avoid confusion with the terminology used in ResNet, we refer to $r^k$ as the right-hand-side of the error equation. In the context of defining a sub-solver, such as within the multigrid algorithm, one may employ several sweeps of iterations within the stationary iterative method.

The Krylov subspace method~\cite{golub2013matrix} represents another iterative framework where the solver $B$ serves as a preconditioner. This method often accelerates the iteration process despite introducing some overheads. Notable examples include
\begin{itemize}
    \item the Conjugate Gradient (CG) method for symmetric positive definite matrices,
    \item the Minimal Residual (MINRES) method for symmetric indefinite matrices,
    \item the Generalized Minimal Residual (GMRES) method for general matrices, and
    \item the Flexible Generalized Minimal Residual (FGMRES) method for non-stationary preconditioners.
\end{itemize}

It's worth noting that our network, designed as a linear operator, can be effortlessly applied into these two iterative frameworks as a solver $B$.

\subsection{Multigrid Algorithm}
\label{sec:multigrid}

The multigrid algorithm, as described in Trottenberg et al.~\cite{trottenberg2000multigrid}, is a category of solvers used for solving linear systems stemming from the discretization of partial differential equations (PDEs). This algorithm is characterized by its multilevel structure, where it performs a series of iterative processes across grids of different resolutions. 

There is efficient exchange of information via two transfer operators: the restriction operator and the prolongation operator. The restriction operator plays a role in transferring the residual from a finer grid to a coarser one, effectively condensing the information. On the other hand, the prolongation operator is responsible for interpolating the solution from a coarser grid back to a finer grid. These dual operations ensures that while each level of the grid performs localized computations, their collective actions contribute to a global effect.

To provide a clearer understanding, a classical two-grid algorithm is exemplified in Algorithm~\ref{alg:twogrid}.
\begin{algorithm}
    \caption{Two-grid method}
    \label{alg:twogrid}
    \begin{algorithmic}[1]
    \State \textit{Setup} coarse grid operator $A_c$ and transfer operators $P$ and $R$.
    \State Compute the residual $r = rhs - A\ sol$ on the finest grid, and solve the correction equation $A\ e = r$ on the finest grid approximately.
    \State Restrict the residual $r$ to a coarser grid $r_c$, i.e., $r_c = R\ r$.
    \State Solve the correction equation $A_c\ e_c = r_c$ on the coarser grid by a coarse grid solver.
    \State Interpolate the correction $e_c$ from the coarser grid to the finest grid, and update the solution $sol$ on the finest grid, i.e., $sol = sol + P\ e_c$.
    \State Do additional approximate correction on the finest grid.
    \end{algorithmic}
\end{algorithm}
The multigrid algorithm extends this method by recursively applying these steps as coarse grid solver on progressively coarser grids. Such extension forms a multigrid V-cycle. To prevent recursive calling, the algorithm can be implemented with a down cycle loop and an up cycle loop. In the context of solving linear equations, enhancing accuracy can be achieved via integrating this multigrid cycle into iterative methods mentioned in Subsection~\ref{sec:itermethods}.

Specifically, we implement a classical geometric multigrid (GMG) method on a structured grid with the following settings:
\begin{itemize}
    \item Reduce the coarse grid size to the half of fine grid in each dimension;
    \item Discretize the PDE on the coarse grid to set up the coarse grid operator $A_c$;
    \item Use the bilinear interpolation to set up the prolongation operator $P$ and restriction operator~$R$;
    \item Use several sweeps of weighed-Jacobi smoother to solve the equation $A\ e = r$ on each level.
\end{itemize}
This GMG method is straightforward to implement and efficient for many elliptic PDEs. However, its efficiency diminishes when applied to more complex problems, such as the convection-diffusion equation. We will use this method as a baseline for comparison with our learnable solver in the experiments.

\subsection{Neural Network}
\label{sec:neuralnetwork}

Deep learning encompasses a diverse array of components with learnable parameters, among which CNNs are particularly notable. CNNs are a category of deep neural networks predominantly used for analyzing visual imagery and are characterized by their shift-invariant or space-invariant properties. In essence, a CNN performs a weighted summation over a local area of the input, applying this same convolutional kernel uniformly across all input regions. The size of the CNN's output can either match or be smaller than the input, which is determined by the stride of the convolution.
Mathematically, this operation can be represented as a matrix $K$ multiplying a vector $x$,
\begin{equation}
    \label{eq:conv}
    y = CNN(x) = K x.
\end{equation}
A transposed operation can be performed via a Transposed Convolutional Neural Network~\cite{dumoulin2016guide} (TCNN) on the input, i.e.,
\begin{equation}
    \label{eq:transconv}
    y = TCNN(x) = K^T x.
\end{equation}
Thus, a TCNN can be used to upscale the input. The combination of CNN and TCNN can emulate the prolongation operator and restriction operator in the multigrid algorithm. We further categorize CNN functionality as follows:
\begin{itemize}
    \item reChannelCNN: It adjusts the input channels to match the output channels while keeping the input size unchanged.
    \item RCNN: Similar to the restriction operator in the multigrid algorithm, it reduces the input size by half, i.e., $stride=2$.
\end{itemize} 
Except for RCNNs and TCNNs, all other CNNs will keep the input size unchanged in our experiments. For simplicity, we use the same symbol for the CNN kernel and its corresponding operator matrix. We use a $3\times 3$ kernel for all CNNs and TCNNs and do not adjust the kernel size in our experiments. These two types of layers form the basic computational units of the proposed solver.

\begin{figure}[htbp]
\centering
\includegraphics[width=0.91\textwidth]{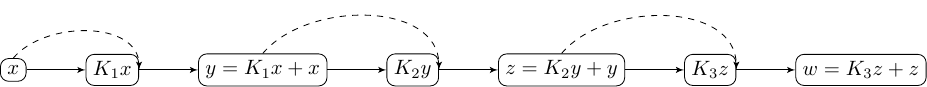}
\caption{Stacked ResNets}
\label{fig:resnet}
\end{figure}

\begin{figure}[htbp]
\centering
\includegraphics[width=0.91\textwidth]{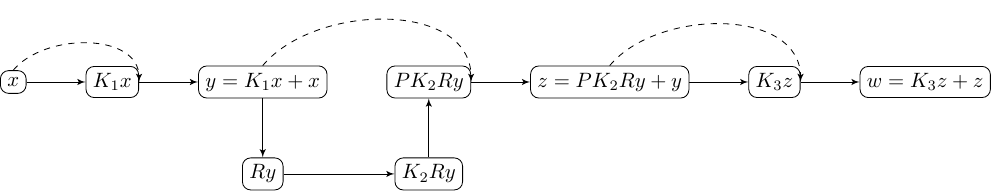}
\caption{Stacked ResNets with twogrid hierarchy}
\label{fig:resnettg}
\end{figure}

An important computation structure in deep learning is ResNet~\cite{he2016deep}, which is a type of CNNs with a skip connection. It was designed to address the vanishing gradient problem in deep neural networks and aligns with numerous state update algorithms in numerical PDEs (e.g., Euler scheme for time-evolution PDEs). Mathematically, it can be expressed as:
\begin{equation}
    \label{eq:resnet}
    y = ResNet(x) = K x + x.
\end{equation}
Such a layer can be stacked to form a deep neural network in eq.~\eqref{eq:deepresnet}, as illustrated in Fig.~\ref{fig:resnet}. This skip connection concept has been successfully used in the related works in their forms.
\begin{equation}
    \label{eq:deepresnet}
    \mathbf{w} = ResNet_3(ResNet_2(ResNet_1(x)))
\end{equation}
With restrict and prolongation operators, we can naturally form a ResNet network with a twogrid hierarchy (see Fig.~\ref{fig:resnettg}). 
To extend the down and up cycle deeper we gain the multigrid structure of the proposed solver (see Fig.~\ref{fig:solve},~Alg.~\ref{alg:solve}).

\begin{remark}[\textbf{State update strategy}]
    In the process of solving a linear system within a stationary iterative method (refer to eq.~\eqref{eq:stationary_iteration}), the solution update differs slightly from the state update of ResNet (see eq.~\eqref{eq:resnet}). This difference lies in the participation of $rhs$ to compute the equations' residual. In the early stage of this research, we developed a multigrid-hierarchical extension that closely aligns with the stationary iterative method. However, this approach did not significantly improve the results and introduced additional convolutional layers. Therefore, we opted for the simpler ResNet as the basis for the proposed solver structure.
\end{remark}

Despite these important networks and structures, another critical component for deep learning is the nonlinear activation function, which combines with the layered structure to form a nice universal approximator. It is a common practice to set nonlinear activation function $\mathbf{\sigma}$ after each CNN layer and add bias term $b$. The formula expression of the nonlinear ResNet is
\begin{equation}
    \label{eq:resnet_nonlinear}
    y = NonLinearResNet(x) = \mathbf{\sigma} (K x + b) + x.
\end{equation}
We will utilize it where nonlinearity is needed. 

\subsection{Discretization of PDEs on Structured Grids}
\label{sec:discretization}
CNNs can implement discretization of a linear PDE in a structured grid. Consider a two-dimensional convection-diffusion equation with heterogeneous diffusion coefficient and a fixed velocity $\mathbf{v}$
\begin{equation}
    \label{eq:convection_diffusion}
    -\mu(x,y) (\partial_{xx} u(x,y) + \partial_{yy} u(x,y)) + \mathbf{v} \cdot \nabla u(x,y) = f(x,y),
\end{equation}
As velocity is set as a constant vector, we use a normalized vector in this paper. The first-order upwind scheme~\cite{larsson2003partial} for this equation can be implemented with two CNNs, as follows:
\begin{equation}
    \label{eq:convection_diffusion_kernel}
    K = \frac{1}{h} \Big\{ 1/\mathbf{Re}
    \begin{bmatrix}
    \ \ 0 &    -1 & \ \ 0 \\
       -1 & \ \ 4 &    -1 \\
    \ \ 0 &    -1 & \ \ 0
    \end{bmatrix} + \mathbf{v}\cdot K_{upwind} \Big\},
\end{equation}
where $h$ is the grid length, $\mathbf{Re} = h/\mathbf{\mu}$ 
is the mesh Reynolds number~\cite{cheng1978computational}, and $K_{upwind}$ is the upwind scheme kernel. We use a bold symbol for $\mathbf{Re}$ to distinguish with its range limit (see Subsection~\ref{sec:basic_settings}). $\mathbf{Re}$ measures the dominance of convection over diffusion, or the ratio of the unsymmetric part from an algebraic perspective. If we keep $\mathbf{Re}$ unchanged for a larger grid, we are considering an equation with lower diffusion coefficient. This equation is often considered difficult to handle if $\mathbf{Re}$ is significantly high~\cite{zhang1997accelerated, gupta1997compact}. $K_{upwind}$ is a weighted summation of several space-difference kernels below,
\begin{equation}
    dx^+ = \begin{bmatrix}
    0 & \ \ 0 & 0 \\
    0 &    -1 & 1 \\
    0 & \ \ 0 & 0
    \end{bmatrix},\ 
    dx^- = \begin{bmatrix}
    \ \ 0 & 0 & 0 \\
       -1 & 1 & 0 \\
    \ \ 0 & 0 & 0
    \end{bmatrix},\ 
    dy^+ = \begin{bmatrix}
    0 & \ \ 1 & 0 \\
    0 &    -1 & 0 \\
    0 & \ \ 0 & 0
    \end{bmatrix},\ 
    dy^- = \begin{bmatrix}
    0 & \ \ 0 & 0 \\
    0 & \ \ 1 & 0 \\
    0 &    -1 & 0
    \end{bmatrix}
\end{equation}
Given velocity $\mathbf{v} = (v_x, v_y)$, we can define 
\begin{equation}
    v_x^+ = max(v_x, 0),\ v_x^- = min(v_x, 0),\ v_y^+ = max(v_y, 0),\ v_y^- = min(v_y, 0)
\end{equation}
and then the upwind scheme kernel can be written as
\begin{equation}
    K_{upwind} = \big( v_x^+ dx^- + v_x^- dx^+,\ v_y^+ dy^- + v_y^- dy^+ \big)
\end{equation}

In this study, we only consider problems with zero Dirichlet boundary conditions on a two-dimensional rectangular domain. Correspondingly, we apply zero padding to CNNs when necessary on solve phase. Different padding strategies can be employed for various boundary conditions. It's worth noting that through discretization on the rectangular domain with a structured grid, each variable can form a tensor. During the learning process, we normalize the input by eliminating the factor $1/h$ in kernel $K$ for different problem sizes and refer to $coef = 1/\mathbf{Re}=\mu/h$ as the coefficient~(tensor), which determines the matrix $A$ in the linear system.

\begin{remark}[\textbf{Non-divergence form vs divergence form}]
    \label{rmk:tested_pdes}
    The conventional convection-diffusion equation is typically expressed in a divergence form as follows:
    \begin{equation}
        \label{eq:convection_diffusion_div}
        -\nabla \cdot (\mu(x,y) \nabla u(x,y)) + \mathbf{v} \cdot \nabla u(x,y) = f(x,y).
    \end{equation}
    However, for the sake of simplicity in implementation with CNNs, we use the non-divergence form~\cite{safonov2010non} in equation~\eqref{eq:convection_diffusion}. This choice facilitates efficient $A\ sol$ operations during the solve phase, e.g. within stationary iterative methods, as outlined in Subsection~\ref{sec:itermethods}. It is worth noting that the $A\ sol$ operation will be frequently utilized, especially in the context of the GMG method, whereas its usage in our solver is comparatively much less, owing to the network's independence of equation formula and fast convergence. To ensure a time-efficient implementation of the GMG method, we adopt the non-divergence form.

    The focus of our paper is to develop a general-purpose solver rather than targeting a specific problem. In practice, users have the flexibility to implement their own high-performance kernels to execute $A\ sol$ in iterative methods, which is not limited to CNN implementation.
\end{remark}

\section{Methodology}
\label{sec:methodology}
Keeping the preliminaries in Subsections~\ref{sec:multigrid} and~\ref{sec:neuralnetwork} in mind, in the upcoming subsections, we'll first appreciate the raised three principles and outline the two phases of the network solver. Then we delve into the details of these two phases, providing explanations of their design and demonstrating how they align with our guiding principles. Additionally, we'll discuss the training aspects of our solver, providing insights into how the network is prepared to tackle problems effectively.

\subsection{Design Principles and Network Phases}
\label{sec:designprinciples}
Here, we reiterate the three guiding principles that underpin the design of our proposed solver:
\begin{itemize}
\item \textbf{P1.} Implement a multigrid hierarchical computation architecture for state updates\label{P:1};
\item \textbf{P2.} Design the network to function linearly for the right-hand-side term of the equations\label{P:2};
\item \textbf{P3.} Share weights across different levels within the multigrid hierarchy\label{P:3}.
\end{itemize}
These principles not only resonate with many classical multigrid algorithms but also serve as the foundation for our development of a learnable solver with exceptional functionality, which we will elaborate on soon. The first principle is a common practice for ensuring good convergence in large problems. We then come to analyze the other two principles.

\vspace{0.5em}
\noindent\textbf{P2. Linearity with respect to the right-hand-side.}
The linear nature of the solver network, in its interaction with the right-hand-side (RHS) term, is a pivotal property for both the solve phase and the training process. On one hand, the efficacy of deep learning is closely linked to the distribution of input data. On the other hand, achieving higher accuracy often involves an iterative method, leading to an updated RHS within a possibly new distribution. This necessitates the generation of an appropriate RHS distribution during the iterative solve phase, as demonstrated in prior research~\cite{azulay2022multigrid, lerer2023multigrid}. Consequently, these studies require a classical solver with specialized expertise for such data generation. 

In our methodology, we exploit the linear property of the network to enhance its generalization ability. It is crucial because training data, when linearly combined, can span a much broader data space. Utilizing a white noise dataset for training proves to be adequate due to this property. Furthermore, linearity enables its usage in various iterative methods, extending beyond FGMRES utilized in work~\cite{azulay2022multigrid, lerer2023multigrid}, to improve solution accuracy or accelerate convergence. This strategy not only simplifies the training process to be free of a ready classical solver to generate RHS dataset, but also broadens the applicability and improves the potential robustness of the trained solver.

\vspace{0.5em}
\noindent\textbf{P3. Weights sharing across different levels.}
We aim for our solver to be adaptable to a wide range of problem sizes. To accomplish this goal, we implement weight sharing across different levels of the multigrid hierarchy. This approach ensures that the network can use any number of levels to accommodate various problem sizes. Furthermore, weight sharing facilitates the training process by allowing the network to learn the mapping from small problems and subsequently generalize to larger ones. Therefore, this principle is crucial for enhancing the scalability of the network.

\begin{remark}[\textbf{No weight sharing when number of levels fixed}]
    When dealing with certain special problems, such as the Helmholtz equation, increasing the number of levels in the multigrid hierarchy may not always result in improved performance. A common approach in classical methods is to apply only a few levels of the multigrid hierarchy and use a direct solver for the coarsest grid. This approach has been successfully employed in the work of Azulay et al.\cite{azulay2022multigrid} and Lerer et al.\cite{lerer2023multigrid}. In such cases, the weight-sharing mechanism might not be necessary, and the network remains fixed for different problem sizes.
\end{remark}

\vspace{0.5em}
\noindent\textbf{Network phases.}
Our proposed solver, inspired by the multigrid algorithm (see Alg~\ref{alg:twogrid}), operates in two similar phases: the setup phase and the solve phase. 

\begin{enumerate}
\item Setup Phase: This initial phase involves constructing the necessary information for each level of the multigrid hierarchy. This information lays the foundation for the subsequent solve phase. The setup is performed only once for a specific problem.
\item Solve Phase: In this phase, we utilize the information prepared during the setup phase to sequentially update the state across different levels of the hierarchy. The solve phase is iterative, allowing for repeated execution until the solution reaches satisfactory accuracy.
\end{enumerate}

\subsection{Setup Phase}
\label{sec:setup}
\begin{figure}[htbp]
    \centering
    \includegraphics[width=0.91\textwidth]{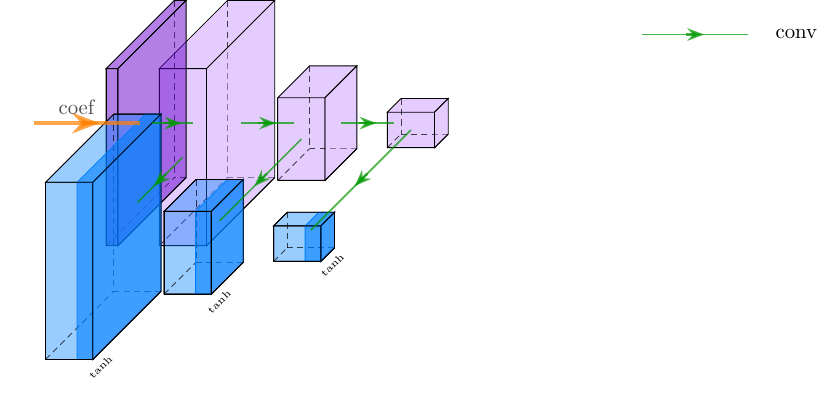}
    \vspace*{-20pt}
    \caption{The 3-level setup network. \footnotesize The inner network (in pink) is used to restrict the problem coefficient tensor (in magenta) from fine to coarse grid. The outer network (in blue) maps these coefficient on each level into setup tensors needed for the solve phase. The right-banded blue boxes indicate the use of the $tanh$ activation function.}
    \label{fig:setup}
\end{figure}

For this structured grid problem, we consider a similar strategy as the GMG method (see Subsection~\ref{sec:multigrid}) to setup the solver, which is to discretize the PDE on the coarser grid to form the coarse grid operator $A_c$. Our setup phase can be split into two steps:
\begin{itemize}
    \item \textbf{Setup1.} use RCNN to restrict the problem coefficient tensor from fine grid to coarse grid
        \begin{equation}
            \label{eq:setup1}
            coef_{l+1} = RCNN(coef_{l}),
        \end{equation}
    \item \textbf{Setup2.} use multi-layered nonlinear ResNet to map these coefficient on each level into setup tensors needed for the solve phase
        \begin{equation}
            \label{eq:setup2}
            setup\_out_{l} = NonLinearResNet(coef_{l}).
        \end{equation}
\end{itemize}
The setup tensors, once computed, are stored and utilized in every iteration of the solve phase. Given that the setup phase is executed only once, the network in \textbf{Setup2} can be designed with sufficient layers to learn complex nonlinear mappings, and we choose four layers in use. Our experiments indicate that the $tanh$ activation function outperforms $ReLU$ for the chosen discretized problem. Consequently, we use $tanh$ after every CNN hidden layer in the ResNet of \textbf{Setup2}, the sole part where we aim to construct nonlinear mappings. Fig.~\ref{fig:setup} illustrates the setup computation framework, and Alg.~\ref{alg:setup} provides a detailed implementation.
\begin{algorithm}
    \caption{Setup phase}
    \label{alg:setup}
    \begin{algorithmic}[1]
    \State \textbf{Input:} $coef$ is the coefficient tensor, $level$ is the number of levels in the multigrid hierarchy.
    \State \textbf{Output:} $setup\_outs$ is the list of setup tensors for each level.
    \Procedure{Setup}{$coef, level$}
        \State $coef_1 \gets reChannelCNN(coef)$ \Comment{Setup1. coefficients on finest grid}
        \For{$l$ in $1, 2, ..., level$}
            \State $setup\_out_l \gets NonLinearResNet(coef_l)$ \Comment{Setup2. map coefficients to setup tensors}
            \If{$l < level$}
                \State $coef_{l+1} \gets RCNN(coef_l)$ \Comment{Setup1. restrict coefficients to coarser grid}
            \EndIf
        \EndFor
        \State \Return $setup\_outs \gets [setup\_out_1, setup\_out_2, ..., setup\_out_{level}]$
    \EndProcedure
    \end{algorithmic}
\end{algorithm}

To use the same model for different levels, we share the corresponding weights of CNNs for the two multigrid hierarchies, i.e., following \textbf{P3}. An effect of weights sharing is the same number of channels over various levels.

\begin{remark}[\textbf{Matrix coefficients as input}]
    Not only problem coefficient but any tensor containing complete information used to define the linear system operator $A$ can be imported into the setup network. Therefore, when the discretization method is unknown and only matrix data is accessible, we can import matrix data as a tensor to the setup network. For example, if the matrix dataset is of the stencil form
    \begin{equation}
        \label{eq:kernel}
        K =
        \begin{bmatrix}
        0 & 0 & 0 \\
        0 & c_3 & c_2 \\
        0 & c_1 & 0
        \end{bmatrix},
    \end{equation}
    where $c_1, c_2, c_3$ are probably different over the domain, we can import this matrix as coefficient tensor
    \begin{equation}
        \label{eq:tensor}
        concat(c_3, c_2, c_1),
    \end{equation}
    to the setup network. This way, we can use the proposed setup network to learn the mapping from the matrix to setup tensors.
\end{remark}

\begin{remark}[\textbf{Setup phase structure}]
    There might be lots of structure for the setup phase to explore without conflict with the three principles. Work~\cite{azulay2022multigrid} use a multigrid-like V-cycle and work~\cite{lerer2023multigrid} goes through a down cycle. We have tried several and finally chose this one to mimic the process of the GMG setup phase. It works slightly better than others. Furthermore, we can see that there is no computation dependence on \textbf{Setup2} network between different levels, which is of potential high complexity. Therefore one can optimize code to further parallelize \textbf{Setup2}. Though there is much flexibility, it works as the handler of the system coefficient data and might be closely related to the final learning performance. We will leave exploration of the impact of the setup phase structure to future work.
\end{remark}

\subsection{Solve Phase}
\label{sec:solve_phase}
\begin{figure}[htbp]
    \centering
    \includegraphics[width=0.91\textwidth]{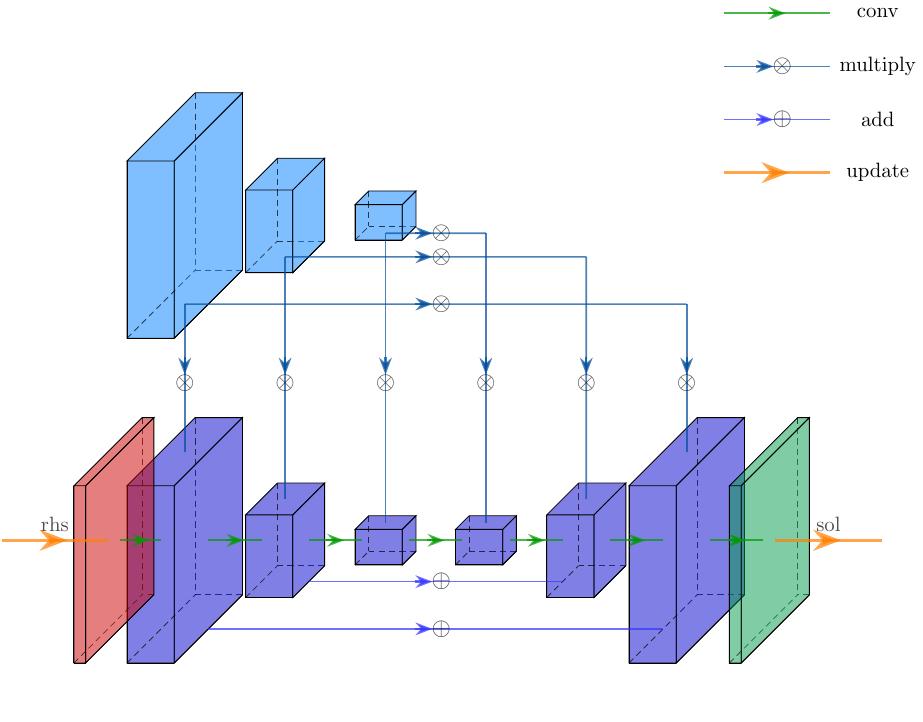}
    \vspace*{-20pt}
    \caption{The 3-level solve network. \footnotesize The top part displays the setup output tensors (in blue). The bottom (in purple) represents the solve phase, where the right-hand-side and the solver output are updated within a user-defined iterative method.}
    \label{fig:solve}
\end{figure}

Solve phase is the central part of our design, reflecting all three principles. To comply with \textbf{P1}, readers should have a comprehensive understanding of the features of different neural networks and the architecture of the two-grid ResNet (see Fig.~\ref{fig:resnettg}). Our solve phase directly extends the two-grid ResNet architecture in a multigrid-hierarchical fashion to update the state. The specific implementation is detailed in Alg.~\ref{alg:solve}, with a down cycle loop to sweep from the finest to coarsest level and a up cycle loop back to finest level. A 3-level network illustration is in the lower part of Fig.~\ref{fig:solve}, with the setup output tensors displayed at the top, which we will discuss soon.
\begin{algorithm}[htbp]
\caption{Solve phase}
\label{alg:solve}
\begin{algorithmic}[1]
\State \textbf{Input:} $setup\_outs$ is the list of setup tensors for each level, $rhs$ is the right-hand-side tensor, $level$ is the number of levels in the multigrid hierarchy.
\State \textbf{Output:} $sol$ is the solution tensor.
\Procedure{Solve}{$setup\_outs, rhs, level$}
\State $x_1 \gets reChannelCNN(rhs)$
\For{$l$ in $1, 2, ..., level$} \Comment{Down cycle loop}
    \State $x_l \gets ResNetDown(setup\_outs[l], x_l)$ \Comment{Update state on each level}
    \If{$l < level$}
        \State $x_{l+1} \gets RCNN(x_l)$ \Comment{Restrict state to coarser grid}
    \EndIf
\EndFor
\For{$l$ in $level, level-1, ..., 1$} \Comment{Up cycle loop}
    \State $x_l \gets ResNetUp(setup\_outs[l], x_l)$ \Comment{Update state on each level}
    \If{$l > 1$}
        \State $x_{l-1} \gets x_{l-1}+TCNN(x_l)$ \Comment{Add state from coarser grid}
    \EndIf
\EndFor
\State $sol \gets reChannelCNN(x_1)$
\State \textbf{return} $sol$
\EndProcedure
\end{algorithmic}
\end{algorithm}

To adhere to \textbf{P2}, we avoid using any nonlinear activation function or bias in the solve phase. The trick lies in leveraging the setup tensors while maintaining the network as a linear operator for the right-hand-side. We find that the simple \textbf{multiplication} operation works well. This means that the input of a convolutional layer in $ResNetDown()$ and $ResNetUp()$ is multiplied with a setup tensor on each level, i.e.,
\begin{align} 
    \label{eq:heteroresnet}
    y &= ResNetDown(setup\_out, x) = K_{down} (setup\_out \times x) + x, \\
    y &= ResNetUp(setup\_out, x)   = K_{up} (setup\_out \times x) + x,
\end{align}
where $K_{down}$ and $K_{up}$ are the convolutional kernels in $ResNetDown()$ and $ResNetUp()$, respectively, and $\times$ denotes the element-wise multiplication. 
From a GMG perspective, we construct a smoother for each level based on its problem coefficient. Lastly, \textbf{P3} is easily obeyed by sharing weights. In other words, the functions $ResNetUp()$, $RCNN()$, $TCNN()$, and $ResNetDown()$ are the same across different levels. 

The right-hand-side and the solver output are updated within a user-defined iterative method. These operations above form the entire landscape of the solve phase, and an illustration is in Fig.~\ref{fig:solve}. 
Like the multigrid method, where the user can define the smoother \textit{sweeps} on each level, we can also specify the number of layers in $ResNetDown()$ and $ResNetUp()$. To minimize the serial computation of stacked ResNets, we use just one sweep on each level in our experiments.

\FloatBarrier
We now provide a series of titled remarks to discuss numerous additional insights, which are beneficial for a deeper understanding of the proposed solver and can inspire further research. We encourage readers to select their most interested ones to explore.

\begin{remark}[\textbf{Black-box solver}]
    While we refer to our solver as a ``learnable multigrid solver'', it's important to note that only the computational structure and multilevel information transfer pattern bear similarity to classical multigrid solvers. Unlike classical multigrid solvers, our solver lacks explicit mathematical explanations for individual components. Instead, all components collectively form a black-box solver for the user, and training through data to achieve small loss can ensure good performance.
\end{remark}

\begin{remark}[\textbf{Symmetric problem}]
    In some cases, the problem may have an inherent algebraic structure that the user wishes to preserve when setting up the solver. For example, for a symmetric problem, one could design a solver that operates as a symmetric operator and use the MINRES Krylov method to accelerate convergence. In our framework, we can enforce the network to be symmetric by sharing weights between $ResNetDown()$ and $ResNetUp()$, and between $RCNN()$ and $TCNN()$.
\end{remark}

\begin{remark}[\textbf{Non-heterogeneous problem}]
    We can simplify both the setup and solve phases for non-heterogeneous problems, where the problem coefficient remains constant. In these scenarios, we only need to generalize the solver over a range of one or several parameters. The space-invariant property of CNNs allows for direct use of the CNN layer as the linear operation, eliminating the need for our proposed multiplication operation with a setup tensor. Consequently, we can use a simple fully-connected neural network instead of a $NonLinearResNet()$ to set up the problem information for solve phase. The output of the setup phase then becomes a column vector for each level, which can be set as the kernel weights of the convolutional layers in $ResNetDown()$ and $ResNetUp()$. This approach coincides with the concept of meta-learning, a multi-task learning framework, as used in Meta-MgNet~\cite{CHEN2022110996}.
\end{remark}

\begin{remark}[\textbf{Other multigrid patterns}]
    The proposed structure represents just one of several options that adhere to our three guiding principles. Principle \textbf{P1} establishes the basic computation pattern. We can emulate other multigrid patterns, such as semi-coarsening, or employ additive-type operations instead of multiplicative-type ones. These variations could offer additional benefits, similar to those observed in classical multigrid methods. Exploring these alternative structures for the solve phase remains a subject for future research.
\end{remark}

\begin{remark}[\textbf{Vector-type PDEs}]
    Our proposed solver can be naturally generalized to accommodate vector-type PDEs by simply adjusting the number of input and output channels. However, it does not directly support PDEs where different unknowns vary in size, as exemplified by the Marker-And-Cell (MAC) scheme for Stokes equations~\cite{chen2016finite}. Addressing this limitation and extending the solver's capabilities to handle such problems is a goal for our future research.
\end{remark}

\subsection{Training}
\label{sec:training}
The training of the proposed solver is minimizing the square of the residual norm of the linear equations, whose square root is always taken as the convergence criterion of iterative methods. The loss function can be written as
\begin{equation}
    \label{eq:loss}
    L = \frac{1}{N}\sum_{i=1}^{N} ||rhs_i - A_i sol_i||^2
\end{equation}
where $N$ is the number of data, $rhs_i$ is the input (right-hand-side), $sol_i$ is the output (solver-guessed solution), and $A_i$ is the discretized linear operator, which is decided by another input -- coefficient tensor $coef_i$. Under the machine learning context, it is an unsupervised regression task, as we do not know the ground truth of the solution. Considering the frequency distribution of the white noise process is still white noise, which means it covers all ranges of frequency, we generate $rhs_i$ as a white noise tensor. In practice, the coefficient tensor $coef_i$ should be generated to fit one's application scenarios, and our generating method for experiments is in Subsection~\ref{sec:basic_settings}.

\begin{figure}[htbp]
    \centering
    \includegraphics[width=0.91\textwidth]{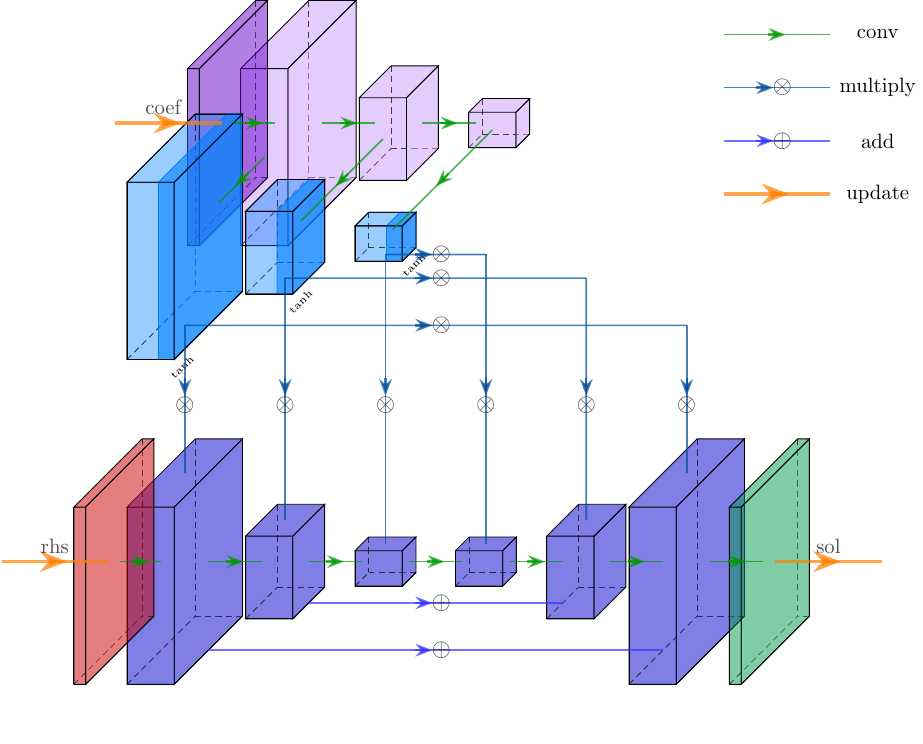}
    \vspace*{-20pt}
    \caption{The 3 level training network. \footnotesize A combination of setup and solve networks. The right-hand-side, coefficient and solver output are updated in every training batch.}
    \label{fig:train}
\end{figure}

During training, we only perform one iteration of the solve phase for each coefficient tensor. The right-hand-side and coefficient are generated randomly for every batch. The whole training process, which combines the setup and solve networks, is illustrated in Fig.~\ref{fig:train}. 

In practice, iterative solvers are typically employed to solve large linear systems. While the training process is assumed to be completely offline, directly training the solver on the large systems can still pose significant challenges in terms of time and memory consumption. Taking the benefit of its applicability to different problem sizes with different levels, a more feasible strategy is to train the solver network on smaller systems with fewer levels. Note that we train solver network on linear systems of varying smaller sizes, after which the trained network can be applied to different systems of larger sizes. We record the training time and memory consumption in our experiments in Appendix~\ref{sec:training_history}.

Essentially, the training process varies the right-hand-side and coefficient tensors, as well as problem sizes, enabling the network to generalize across different problem instances. This makes the trained network as prepared as a traditional solver to tackle new problems within this class of PDEs. Therefore, training is assumed offline and akin to manually designing a solver for a specific class of problems. The training process is efficient and straightforward, requiring no specific preparation of the right-hand-side dataset (see analysis of \textbf{P2.} in Subsection~\ref{sec:designprinciples}) and only a user-defined distribution for the coefficient dataset.

\begin{remark}[\textbf{Unsupervised vs supervised learning strategy}]
    We could also employ a supervised learning pattern through the following steps:
    \begin{enumerate}
        \item Generate a random tensor as the ground truth solution $true\_sol_i$.
        \item Use $A_i true\_sol_i$ as the right-hand-side, i.e., $rhs_i$ sent to the solver.
    \end{enumerate}
    In this case, we can define the loss function as the square of the error norm:
    \begin{equation}
        L = \frac{1}{N}\sum_{i=1}^{N} ||sol_i - true\_sol_i||^2
    \end{equation}
    However, generating a suitable $true\_sol_i$ might be questionable. Mathematically, the right-hand-side function usually has weak smoothness and only needs to be of $\mathbb{L}_2$ regularity, while the solution function is usually of strong smoothness and needs to be of, for instance, $\mathbb{H}^1$ regularity or even $\mathbb{H}^1_0$. Therefore, it is not straightforward to generate a suitable $true\_sol_i$ for a given $A_i$. In our experiments, we find that the unsupervised learning pattern works well and does not gain benefit from the supervised learning pattern.
\end{remark}

\section{Numerical Experiments}
\label{sec:experiments}
In this section, we present a series of experiments to evaluate the performance of our proposed solver. In the upcoming subsections, we will first introduce the basic settings of our experiments in Subsection~\ref{sec:basic_settings}. Then, we will explore the impact of the hyperparameter \textbf{hidden channels} and compare the convergence and time consumption of our solver with GMG in Subsection~\ref{sec:hyperparameter}. Subsequently, our aim is to investigate its dependence on coefficient data distribution in Subsection~\ref{sec:datadistribution}, including the effects of \textbf{distribution range} and \textbf{distribution pattern}.

\subsection{Basic Settings}
\label{sec:basic_settings}
Our experiments were conducted on a V100 GPU with 16GB memory, using PyTorch~\cite{paszke2019pytorch} as our deep learning framework. The code is available at \url{https://gitee.com/xiehuohuo77/heteromgcnn.git} (open after publication). The GMG method, also implemented in PyTorch, is used as a comparison, representing a classical method with a similar computation pattern, utilizing identical hardware and software environments, and under the same settings, free of expertise.

\begin{remark}[\textbf{AMG vs GMG}]
    AMG is another popular multigrid method, while there are some essential differences between AMG and our proposed solver. AMG offers a more flexible coarsening pattern and entails a more complex setup phase, which is typically preferred to conduct on CPU. We have experimented with BoomerAMG in HYPRE~\cite{falgout2002hypre} on the tested linear systems. When run on CPU (Intel(R) Xeon(R) Platinum 9242 CPU @ 2.30GHz) with manually optimized number of processes (20) for a grid size of $1023 \times 1023$, it takes approximately 772 ms in total (182 ms for setup and 590 ms for solve), which exemplifies the efficiency of our implementation of GMG~(see Table~\ref{tab:channels_double}).
\end{remark}

\vspace{0.5em}
\noindent\textbf{Problem settings.} We trained the solver on linear systems from the convection-diffusion equation with heterogeneous diffusion coefficient (see Subsection~\ref{sec:discretization}). The convection-diffusion equation has fixed velocity $v = (sin(0.5), cos(0.5))$ over the domain. In this experiment, we generate $coef$ from a given random tensor to range $[\epsilon, 1], 0 < \epsilon < 1$ via two steps:
\begin{enumerate}
    \item firstly linearly map the random tensor to range $[0, -\log_{10}(\epsilon)]$ as $random\_power$,
    \item then map the tensor to range $[\epsilon, 1]$ via $10^{-random\_power}$.
\end{enumerate}
The random tensor is generated as white noise or from other distributions, which we will introduce to discuss the data distribution dependence issue in the experiment Subsection~\ref{sec:datadistribution}. Recall that $\mathbf{Re} = 1/coef$, therefore the range of $\mathbf{Re}$ is $[1, 1/\epsilon]$. For simplicity, we use the non-bold symbol $Re$ to denote the range limit of $\mathbf{Re}$, i.e., $\mathbf{Re} \in [1, Re]$. Unless otherwise specified, we use $Re=1000$ in our experiments.

\vspace{0.5em}
\noindent\textbf{Iterative method settings.} 
In our experiments, each test result represents the average outcome obtained from conducting the test process ten times. For each test, coefficients are randomly selected, ensuring diverse inputs. The right-hand-side (RHS) for these tests is set as a constant, \textbf{1}, which notably is not part of the training dataset for the RHS. Although it's common to employ the Krylov method, particularly the GMRES method, to accelerate convergence, we found that GMRES only marginally reduces the number of iterations for both our proposed solvers and the GMG method in this specific problem. Moreover, the overheads in GMRES lead to increased total computation time. Given this observation, we decided against using the Krylov method in our experiments. Instead, we opted for the \textit{stationary iterative method} for all tested solvers. For assessing convergence, we set the criterion to be the relative residual norm, with a tolerance threshold of 1E-8 for double precision computations. These settings ensure a consistent and fair comparison across different methods under test.

\vspace{0.5em}
\noindent\textbf{MGCNN settings.} We gradually train the solver from minor to larger sizes, simultaneously reducing the batch sizes. The training process covers grids of sizes $31 \times 31$, $63 \times 63$, $127 \times 127$, $255 \times 255$, and $511 \times 511$. We provide the details of the training settings in Appendix~\ref{sec:training_settings}. We tested all solvers over grids ranging from $31 \times 31$ to $4095 \times 4095$, the maximum size our GPU memory can accommodate.  The main hyperparameter of the proposed solver is $level$, $sweeps$ and $(hidden)\ channels$. We found that the reduction of $level$ significantly deteriorates convergence, and one sweep is enough to reach satisfactory convergence. Therefore, we fix $sweeps=1$ and use $level = 4, 5, 6, ...$ with respect to $grid\_size = 31, 63, 255, ...$ for below experiments of the proposed solver. If not specified, we use $channels = 8$ for the proposed solver.

\begin{remark}[\textbf{Training time}]
    \label{rmk:train_speed}
    The training strategy is time-efficient, taking around 2.0 hours on double precision. Training on single precision can be even faster, taking no more than 20 minutes. However, this results in some loss of convergence performance, especially for large grid sizes. Therefore, we only present the outcomes of the solver trained in double precision here. It's important to emphasize again that the training process is entirely offline, as the solver exhibits the ability to generalize across different right-hand-side terms, PDE coefficients, and grid sizes.
\end{remark}

\vspace{0.5em}
\noindent\textbf{GMG settings.} We manually tune the GMG parameters to improve its results. We set $weight=0.67$, perform \textit{three $sweeps$} for the weighed-Jacobi smoother, and use $level = 2, 3, 4, ...$ corresponding to $grid\_size = 31, 63, 255, ...$ for the comparisons below. This setting of $level$ is smaller than the proposed solvers, as GMG does not gain benefit from more levels as the proposed solvers does.

\subsection{Impact of Channels}
\label{sec:hyperparameter}
In this subsection, our focus is the impact of $channels$ on the solver. Specifically, we select $channels = 4, 8, 12$ and present the iteration number and consumed time results in Table~\ref{tab:channels_double}. It's important to note that the time consumption is measured in \textit{milliseconds} (ms). In Table~\ref{tab:channels_float}, we also conduct tests on single precision, using 1E-4 as the convergence tolerance. However, due to precision limit, both GMG and our proposed solvers can only reach the 1E-3 tolerance for the $4095 \times 4095$ grid, so we only present results up to $2047 \times 2047$ grid for the single precision computation. We can see that all of the proposed solvers perform better than GMG in both precision in the solve phase. However, GMG almost takes no time in the setup phase, while the proposed solver takes a higher cost setup phase. Considering that the setup phase is only executed once, and its time cost is much smaller than the iterative solve phase, it is acceptable.
\begin{table}[htbp]
    \centering
    \caption{Impact of channels on the solver (single precision, convergence tolerance is 1E-4), iters: iteration number, setup: setup time (\textbf{ms}), solve: solve time (\textbf{ms})}
    \label{tab:channels_float}
    {\footnotesize
    \begin{tabular}{|c|l|r|r|l|r|r|l|r|r|l|r|r|}
    \hline
    float & \multicolumn{3}{c|}{channel4} & \multicolumn{3}{c|}{channel8} & \multicolumn{3}{c|}{channel12} & \multicolumn{3}{c|}{GMG} \\ \hline
    grid & iters & setup & solve & iters & setup & solve & iters & setup & solve & iters & setup & solve \\ \hline
        31x31 & 5.7  & 2.6  & 9.6  & 4 & 2.4  & 6.8  & 4 & 2.0  & 5.6  & 10 & 1.9  & 39.4  \\ 
        63x63 & 6 & 2.9  & 10.9  & 4 & 2.6  & 6.8  & 4 & 2.6  & 6.6  & 11.3 & 2.4  & 53.8  \\ 
        127x127 & 6 & 3.1  & 11.3  & 4 & 2.8  & 7.0  & 4 & 3.0  & 7.2  & 13 & 2.7  & 68.0  \\ 
        255x255 & 7 & 3.8  & 15.9  & 4.8 & 3.6  & 10.4  & 4 & 3.8  & 9.0  & 15 & 3.7  & 105.9  \\ 
        511x511 & 9 & 4.4  & 23.6  & 5 & 4.0  & 13.6  & 5 & 4.4  & 15.3  & 18 & 4.5  & 150.4  \\ 
        1023x1023 & 13 & 5.1  & 49.6  & 7 & 5.4  & 32.6  & 5 & 7.3  & 29.7  & 24 & 5.2  & 244.3  \\ 
        2047x2047 & 18 & 10.8  & 168.9  & 9 & 17.9  & 135.2  & 6 & 26.0  & 126.6  & 34 & 6.2  & 555.7  \\ \hline
    \end{tabular}
    }
\end{table}
\begin{table}[htbp]
    \centering
    \caption{Impact of channels on the solver (double precision, convergence tolerance is 1E-8), iters: iteration number, setup: setup time (\textbf{ms}), solve: solve time (\textbf{ms})}
    \label{tab:channels_double}
    {\footnotesize
    \begin{tabular}{|c|l|r|r|l|r|r|l|r|r|l|r|r|}
    \hline
    double & \multicolumn{3}{c|}{channel4} & \multicolumn{3}{c|}{channel8} & \multicolumn{3}{c|}{channel12} & \multicolumn{3}{c|}{GMG} \\ \hline
        grid & iters & setup & solve & iters & setup & solve & iters & setup & solve & iters & setup & solve \\ \hline
        31x31 & 11 & 1.7  & 12.3  & 7 & 2.1  & 9.1  & 7.2  & 2.3  & 10.7  & 16.7  & 1.5  & 50.9  \\ 
        63x63 & 12 & 1.7  & 13.1  & 7.1 & 2.3  & 10.5  & 7.1 & 2.6  & 11.7  & 19.2 & 2.0  & 74.8  \\ 
        127x127 & 12 & 2.0  & 15.2  & 8 & 2.8  & 14.4  & 7 & 3.0  & 14.2  & 21 & 2.7  & 104.5  \\ 
        255x255 & 13 & 3.8  & 20.9  & 8 & 3.5  & 18.9  & 7.5 & 3.8  & 20.1  & 25.9 & 3.7  & 173.8  \\ 
        511x511 & 17 & 2.8  & 42.0  & 9 & 4.1  & 33.6  & 8 & 4.7  & 36.5  & 33.7 & 4.5  & 266.2  \\ 
        1023x1023 & 23 & 6.6  & 152.1  & 12 & 10.8  & 123.0  & 9 & 15.5  & 122.3  & 45 & 5.4  & 533.6  \\ 
        2047x2047 & 30 & 24.4  & 687.0  & 15 & 40.6  & 536.6  & 12 & 58.8  & 592.7  & 61 & 5.9  & 2.1E3  \\ 
        4095x4095 & 43 & 96.6  & 3.6E3 & 22 & 161.6  & 3.1E3 & 15 & 234.5  & 3.2E3 & 84 & 8.0  & 1.1E4 \\ \hline
    \end{tabular}
    }
\end{table}

\begin{figure}[htbp]
    \centering
    \subfloat{\includegraphics[width=0.45\textwidth]{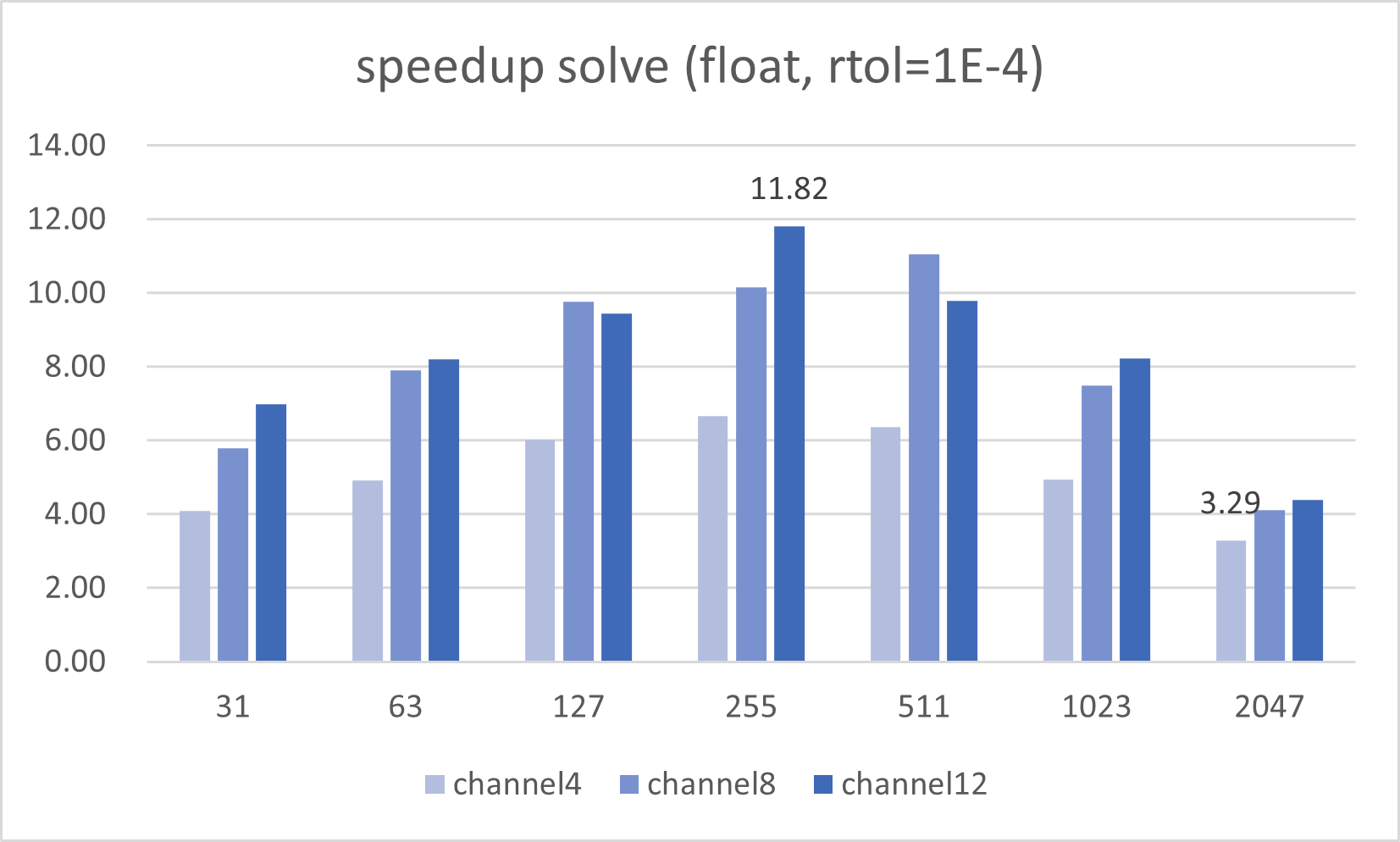}\label{fig:speedup_float}}
    \hfill
    \subfloat{\includegraphics[width=0.45\textwidth]{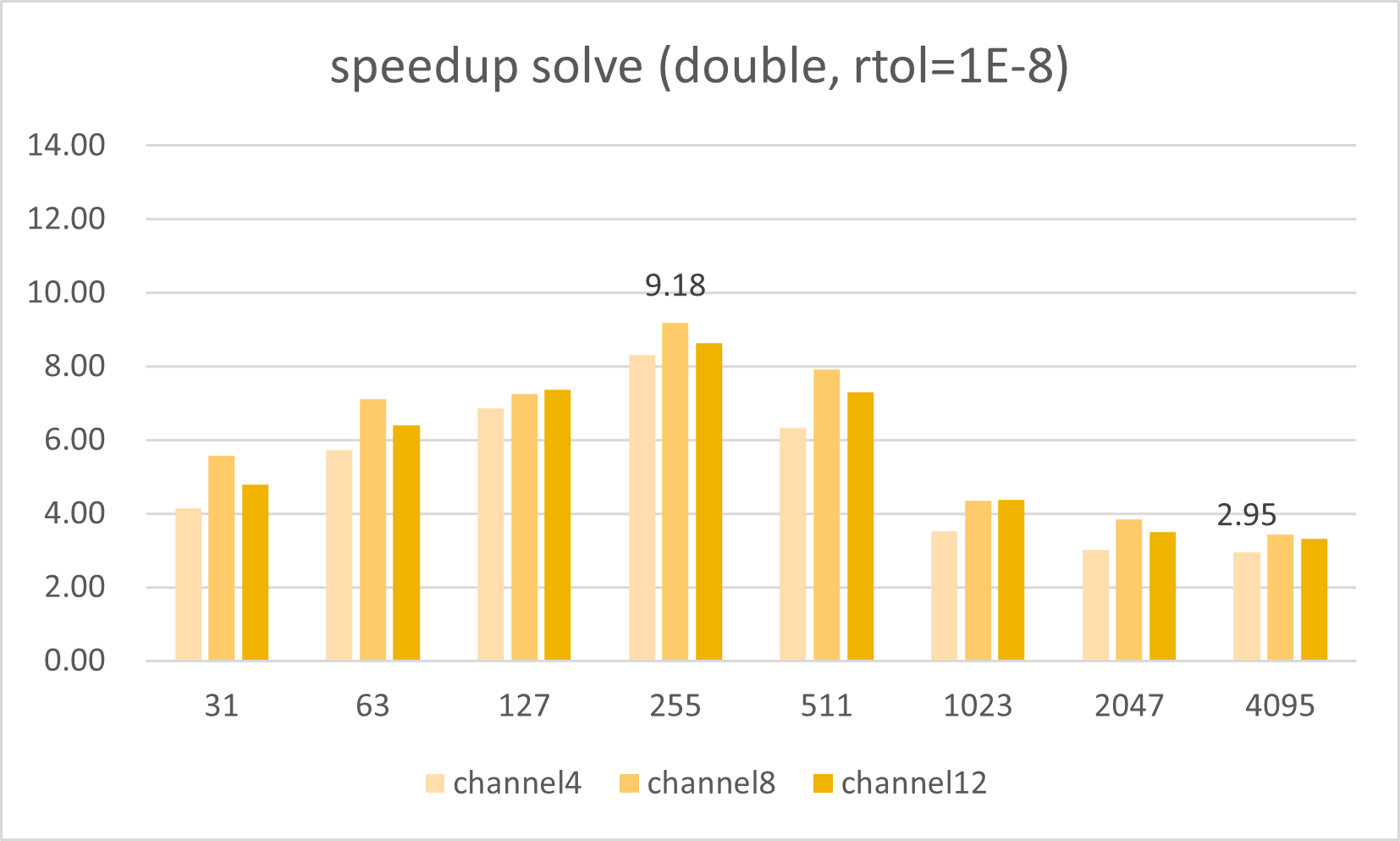}\label{fig:speedup_double}}
    \caption{Time speedup of our trained solver over the GMG method in the solve phase. The left is for single precision, and the right is for double precision.}
    \label{fig:speedup}
\end{figure}
\begin{figure}[htbp]
    \centering
    \subfloat{\includegraphics[width=0.45\textwidth]{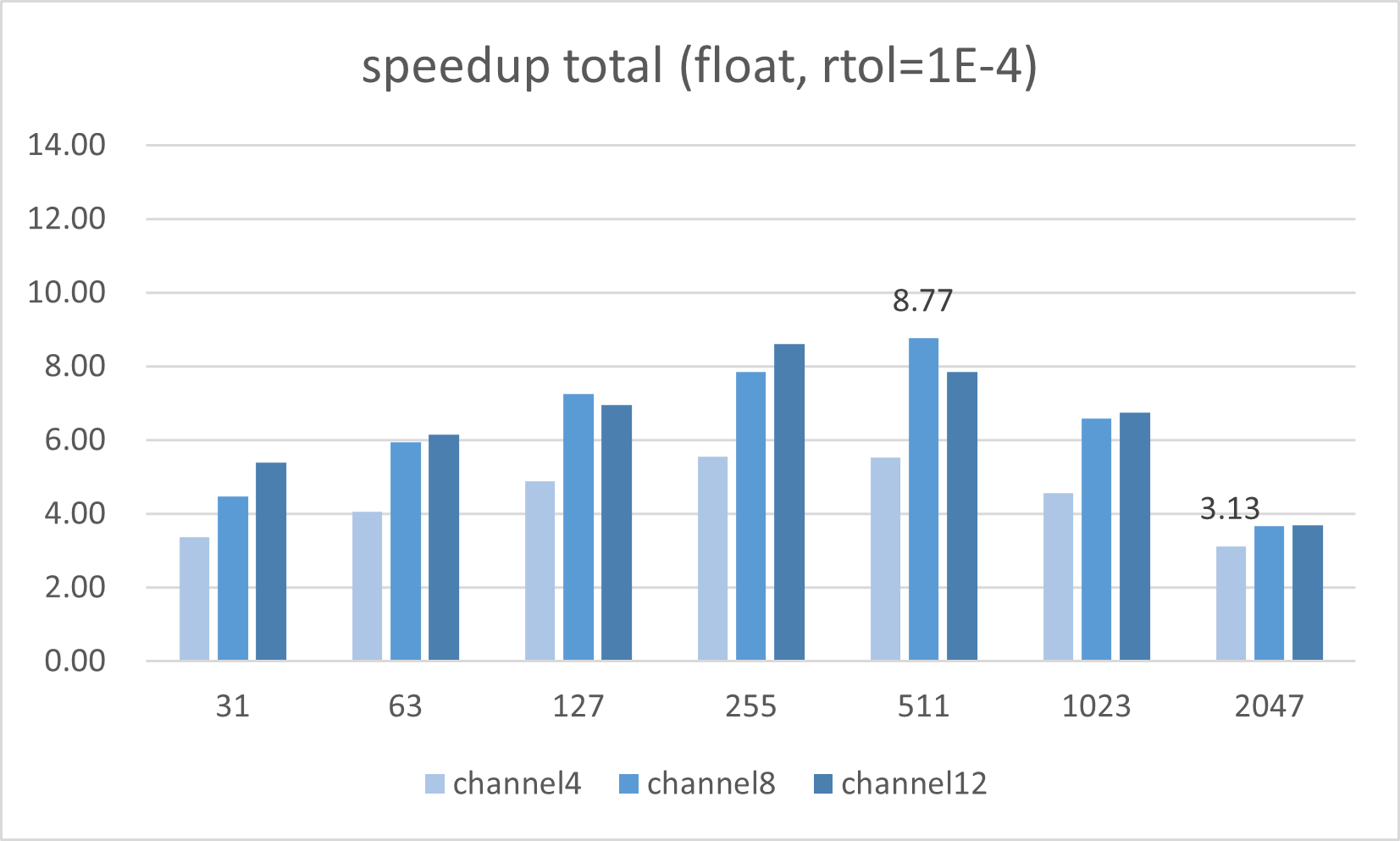}\label{fig:speedup_float_total}}
    \hfill
    \subfloat{\includegraphics[width=0.45\textwidth]{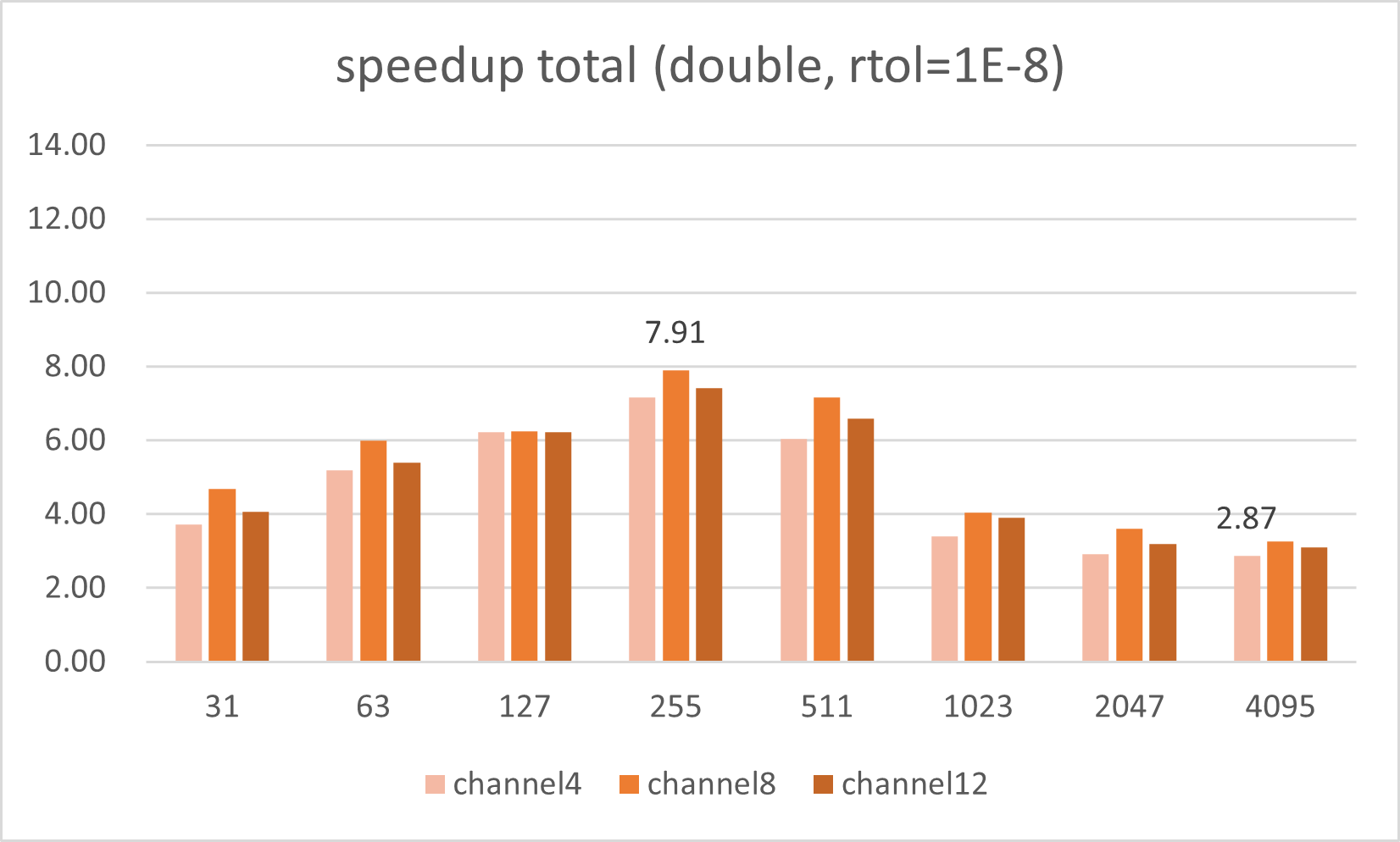}\label{fig:speedup_double_total}}
    \caption{Time speedup of our trained solver over the GMG method in total. The left is for single precision, and the right is for double precision.}
    \label{fig:speedup_total}
\end{figure}

\FloatBarrier
The speedup of the proposed solver over the GMG method in solve phase is shown in Fig.~\ref{fig:speedup}, and the speedup for total time is shown in Fig.~\ref{fig:speedup_total}. In the iterative solve phase, our trained solvers achieve a speedup of around 3x to 10x for the three tested channel numbers on both single and double precision computation. Considering the total time, the speedup is about 3x to 8x, which is also desirable. The highest speedup occurs on medium grid sizes, where the iterations of our trained solver do not increase much, and the GPU computation time is not sensitive to the number of channels. In single precision, the GPU is even less susceptible to such an increase in computation complexity, achieving a higher speedup than double precision. We should remember that GMG involves three sweeps of the smoother (see Subsection~\ref{sec:basic_settings}). Therefore, the speedup of our trained solver mainly stems from the reduction in the number of iterations, owing to the high number of channels and the excellent learning ability of deep neural networks.

\begin{remark}[\textbf{Subspace update view}]
    When examining the update method in the solve phase (see Subsection~\ref{sec:solve_phase}), the main difference from traditional algorithms lies in the number of channels. The proposed solver updates several channels of state simultaneously, which can be viewed as updating a subspace rather than one direction. While multi-channel update may increase the complexity of the solver, such a parallel strategy can be efficiently executed on GPUs. The convergence of such a ``subspace'' method is currently unknown, and we leave it for future work. 

\end{remark}

\FloatBarrier
\subsection{Dependence on Data Distribution of Coefficient}
\label{sec:datadistribution}
Deep learning performance is widely acknowledged to be linked to input data distribution. While our solver performs equally well on test datasets drawn from the same distribution as the training dataset, there remains a concern regarding performance degradation when the tested problems are from a different coefficient distribution. In this subsection, we delve into the influence of both the coefficient distribution range and pattern on the convergence behavior of our trained solver. For simplicity, we borrow the terminology ``transfer learning'' from machine learning to refer to our tests as ``transfer tests''. It's important to note that, unlike traditional transfer learning, no re-tuning is performed during these tests. The results of the GMG solver for comparison are available in Appendix~\ref{sec:GMG_results}.

We generate the coefficient tensor from a random tensor, which comes from a specified distribution (refer to Subsection~\ref{sec:basic_settings}). In addition to white noise, we employ the following distributions for generating the random tensor:
\begin{itemize}
    \item \textbf{CIFAR10}~\cite{krizhevsky2009learning}: A dataset comprising 60,000 $32\times 32$ color images across 10 classes. We utilize the first channel of the images to generate the random tensor.
    \item \textbf{FMNIST}~\cite{xiao2017fashion}: A dataset of 60,000 $28\times 28$ grayscale images representing 10 fashion categories.
    \item \textbf{MNIST}~\cite{lecun1998mnist}: A dataset of 70,000 $28\times 28$ grayscale images depicting the digits 0 to 9.
    \item \textbf{mldata}:  Random tensors generated from a multi-level random data generation method, details of whose algorithm are outlined in Appendix~\ref{sec:mldata}.
\end{itemize}
The first three datasets consist of real-world images, which will be interpolated to the required grid size during experiments. To introduce more varied random pattern data, but not pure white noise, we designed a multi-level data generation method to create \textbf{mldata}, where pixels in a tensor are interconnected in a multilevel pattern. Additionally, we generate the coefficient tensor from white noise with different distribution ranges, namely $Re = 10$ and $Re = 10^5$. 

The proposed solver is exclusively trained on coefficients from white noise with $Re = 1000$, and its convergence is tested on the aforementioned distributions. The number of iterations is summarized in Table~\ref{tab:coef_dist}. We observe that the performance remains reliable for ranges $Re = 10$, $Re = 10^5$, and datasets CIFAR10, \textbf{mldata}. However, it degrades significantly for the other datasets, especially on large grid sizes.

\begin{table}[htbp]
    \centering
    \caption{Iteration number of our trained solver over different coefficient distributions. \footnotesize Re$10$ and Re$10^5$ refer to $Re=10,10^5$ with coefficients generated from white noise. The solver is trained only on coefficients from white noise with $Re = 1000$}
    \label{tab:coef_dist}
    \begin{tabular}{|c|l|l|l|l|l|l|l|}
    \hline
        grid & noise & CIFAR10 & FMNIST & MNIST & mldata & Re$10$ & Re$10^5$ \\ \hline
        31x31 & 7 & 8.2  & 13.6  & 15.1  & 7.7  & 8.1  & 7 \\ 
        63x63 & 7.1 & 9.8 & 17.3 & 18.4 & 8.1 & 9 & 7.8 \\ 
        127x127 & 8 & 11.5 & 19.2 & 19.6 & 9.3 & 9 & 8 \\ 
        255x255 & 8 & 11.1 & 20.6 & 20.9 & 10.5 & 10 & 8 \\ 
        511x511 & 9 & 13.4 & 22.5 & 24 & 11.8 & 12 & 9.9 \\ 
        1023x1023 & 12 & 15.9 & 25.2 & 31 & 13.8 & 15 & 12 \\ 
        2047x2047 & 15 & 19.3 & 39.7 & 45.2 & 16.9 & 21 & 16 \\ 
        4095x4095 & 22 & 25.1 & 62.3 & 105.8 & 23.2 & 31.7 & 22 \\ \hline
    \end{tabular}
\end{table}

\begin{table}[htbp]
    \centering
    \caption{Iteration number of our trained solver over different coefficient distributions. \footnotesize Re$10$ and Re$10^5$ refer to $Re=10,10^5$ with coefficients generated from white noise. The solver is trained on coefficients from a mixture of white noise and CIFAR10.}
    \label{tab:coef_dist_mixed}
    \begin{tabular}{|c|l|l|l|l|l|l|l|}
    \hline
        grid & noise & CIFAR10 & FMNIST & MNIST & mldata & Re$10$ & Re$10^5$ \\ \hline
        31x31 & 7.1  & 7.6  & 9.8  & 10.0  & 7.3  & 7.6  & 7 \\ 
        63x63 & 7.1 & 7.9 & 9.7 & 10 & 7.5 & 7.1 & 7 \\ 
        127x127 & 7.2 & 8 & 9.7 & 9.6 & 8 & 7.6 & 7.6 \\ 
        255x255 & 8 & 8 & 9.5 & 9 & 8 & 8 & 8 \\ 
        511x511 & 10 & 9.6 & 9.6 & 9 & 9.8 & 8 & 10 \\ 
        1023x1023 & 12 & 12.4 & 11.7 & 10.1 & 12.6 & 10 & 13 \\ 
        2047x2047 & 16 & 15.9 & 13.6 & 12.7 & 16.6 & 13 & 17 \\ 
        4095x4095 & 22 & 22.3 & 19.9 & 17.1 & 22.4 & 17 & 23 \\ \hline
    \end{tabular}
\end{table}
To address the distribution transfer issue, we train the solver with a mixture of white noise and CIFAR10 and test its convergence on the same distributions above. The corresponding iteration numbers are listed in Table~\ref{tab:coef_dist_mixed}. We observe that it performs almost equally well on different distribution patterns and range limits. There is only slight increase of the number of iterations on originally trained white noise with $Re = 1000$.

\begin{figure}[htbp]
    \centering
    \includegraphics[width=1.\textwidth]{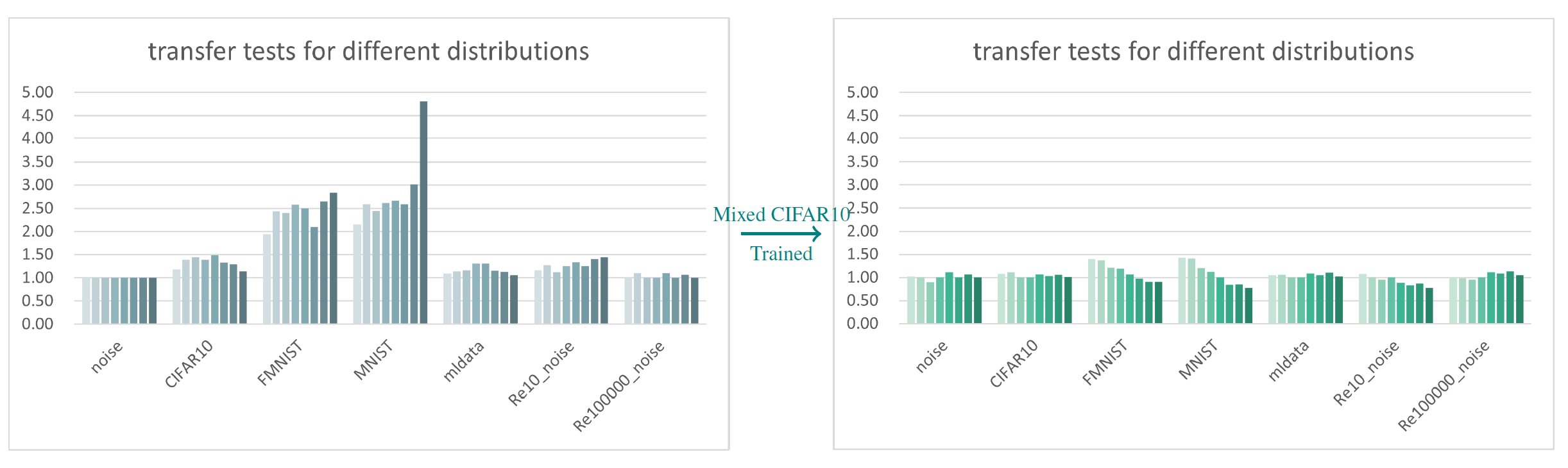}
    \caption{Ratio of iteration numbers of the proposed solver tested on different coefficient distributions over that of white noise. \footnotesize Deeper colors represent larger grid sizes. The left shows results when trained only with white noise, while the right depicts results when trained with a mixture of white noise and CIFAR10.}
    \label{fig:coef_dist_compare}
\end{figure}


To clarify, we refer to the solver trained solely with white noise as the white-noise-trained solver. Conversely, the solver trained on a coefficient dataset from a mixture of white noise and CIFAR10 is termed the mixed-CIFAR10-trained solver. The enhancement achieved by the mixed-CIFAR10-trained solver in transfer tests is clearly depicted in Fig.~\ref{fig:coef_dist_compare}. In this figure, the number of iterations of the solver, when both trained and tested with white noise, is used as a reference point. The ratios of the number of iterations are computed for both the white-noise-trained solver (left) and the mixed-CIFAR10-trained solver (right), each tested on a variety of data distributions. All ratios are reduced to around or even below 1.0, indicating enhanced robustness. For additional detailed results, including the $\log_{10}$ of the coefficient tensors and convergence history of the mixed-CIFAR10-trained solver, refer to Appendix~\ref{sec:our_solver_results}.

In summary, while our proposed solver's performance does show some dependency on the coefficient data distribution, it demonstrates reliable generalization across several untrained distributions. The robustness of the solver can be further improved by training on a mixture of some distributions, with minimal impact on convergence for the originally trained distribution. This suggests that practitioners could benefit from using a broader or mixed distribution of coefficients during training. 

\begin{remark}[\textbf{Training data generation}]
    \label{rmk:train_data}
    In practical terms, we propose the generation of datasets akin to ours, which map image datasets or random data from a generation function to the coefficient tensor, instead of solely relying on real-world coefficient data collection. This method offers convenience and could potentially enhance robustness against unseen coefficient distributions when trained on a mixture of these distributions. One can train the solver on such generated datasets quickly with single precision, and validate its performance on several coefficient data from real applications.
\end{remark}

\section{Conclusion and Future Work}\label{sec:conclusion}

This research presents an innovative deep-learning-based solver designed for sparse linear systems from partial differential equations (PDEs) discretized on structured grids. The solver is built on three key principles: (1) employing a multigrid computational structure, (2) applying the solver linearly to the right-hand-side, and (3) sharing weights across different levels. These principles collectively ensure fast convergence and broad applicability. The multigrid framework enhances the solver's speed and adaptability to large problems. Its linear operation allows for seamless integration into various iterative methods and simplifies the training process, without requiring specialized domain knowledge. Weight sharing across levels enables the solver to adapt to different problem sizes and makes offline training more efficient in terms of both time and memory.

In our study, the solver effectively handles a convection-diffusion equation with varying diffusion coefficients over the domain. We explore the influence of hidden channels within the solver, noting their positive impact on convergence. Due to the GPU's insensitivity to increased channels, particularly in small to medium problem sizes, our solver achieves about 3 to 8 faster performance than the GMG method, under stringent accuracy demands and up to grid size $4095\times 4095$. Furthermore, the transfer tests recommend training the solver on an enlarged distribution of coefficients to enhance its robustness across unseen distributions.

This work marks our initial foray into leveraging AI for solving linear systems from PDEs. Its straightforward deployment and rapid convergence position it as a promising tool. Moreover, it highlights the potential to replace the need for manually designing problem-specific solvers through purely training learnable solvers via data.
Despite the progress detailed in previous sections, further research is necessary. Future directions include extending the principles to sparse linear systems on unstructured grids, developing AI-based nonlinear solvers for nonlinear equations, and addressing the inherent space-invariance of CNNs, which is suboptimal for heterogeneous problems and may lead to time-inefficient computation. We eagerly anticipate continuing our work in this dynamic and evolving field.

\section*{Acknowledgments}
This research is partially supported partly by National Key R\&D Program of China 2020YFA0711900, 2020YFA0711904 and the Strategic Priority Research Program of the Chinese Academy of Sciences, Grant No. XDB0640000.

\appendix
\section{Training}
\label{apd:training}

\subsection{Training Settings}
\label{sec:training_settings}
Although training in this paper does not rigorously follow the dataset concept, we still use the terminology of epoch to control the training process, and the network will go through 1000 batches for each epoch. Our approach involves progressively training the solver on grids of increasing sizes, going through $31 \times 31$, $63 \times 63$, $127 \times 127$, $255 \times 255$, and $511 \times 511$ grids. The maximum batch size is 16 for the smallest size $31 \times 31$ and reduced to half for every doubled size until the minimum batch size 2. The training parameter settings are shown in Table~\ref{tab:training_settings}.
\begin{table}[htbp]
    \centering
    \caption{Training settings}
    \label{tab:training_settings}
    \begin{tabular}{|c|c|c|}
        \hline
        \textbf{Parameter} & \textbf{Value} & \textbf{Description}\\
        \hline
        epochs & 50 & -- \\
        num & 1000 & number of batches in one epoch \\ \hline
        lr & 0.003 & initial learning rate \\
        optimizer & Adam & step\_size=2, gamma=0.8 \\ \hline
        size\_step &  10 & change data\_size, level, batch\_size every size\_step epochs\\
        size &  31 & initial data size, double sizes every size\_step \\
        level &  4 & initial level, increase 1 every size\_step \\
        batch\_size & 16 & initial batch size, reduce to its half every size\_step \\
        max\_size &  511 & maximum data size\\
        min\_batch\_size & 2 &  maximum batch size\\
        \hline
    \end{tabular}
\end{table}

\subsection{Training History}
\label{sec:training_history}
\begin{figure}[htbp]
    \centering
    \subfloat[Loss]{\includegraphics[width=0.45\textwidth]{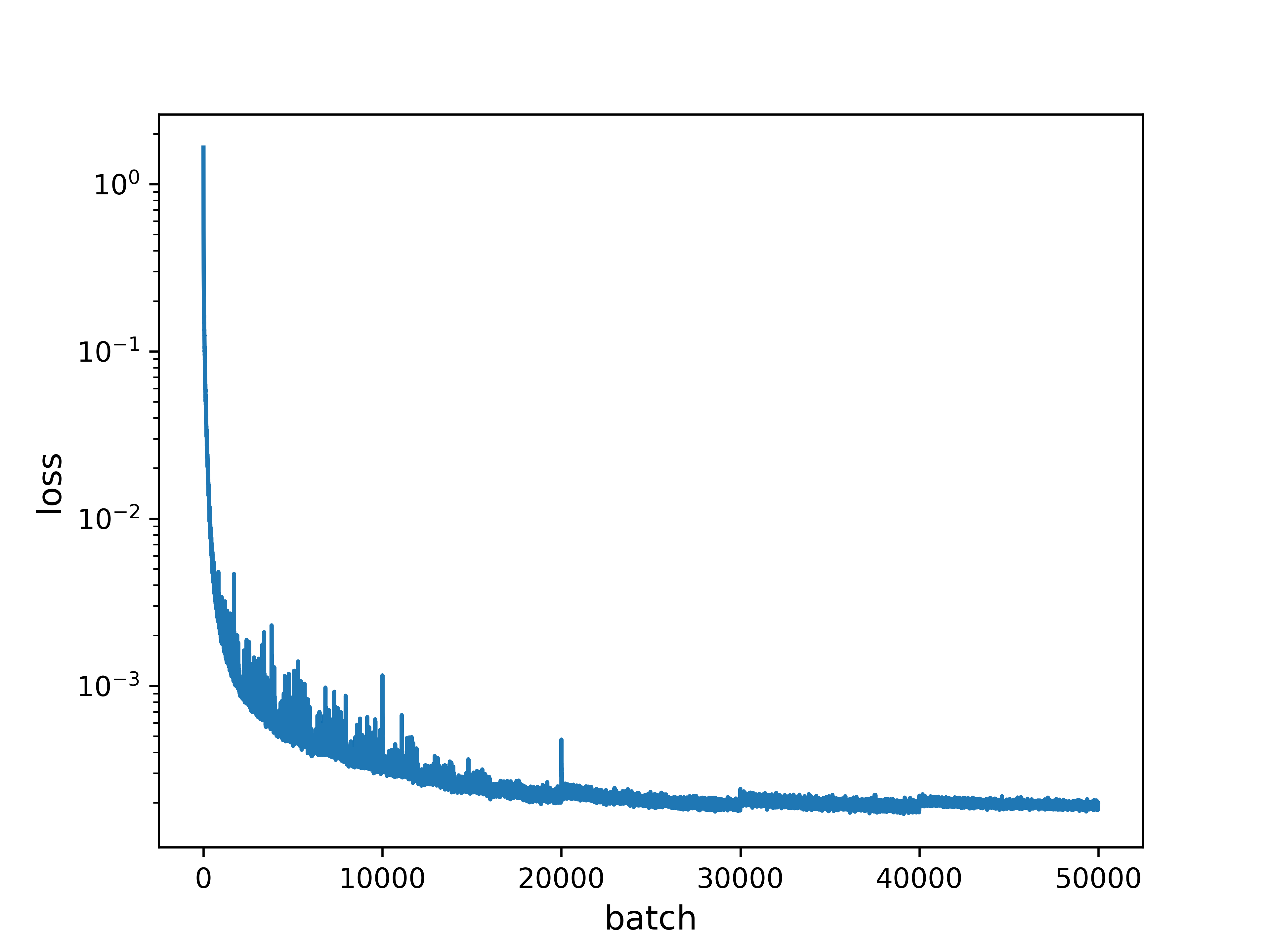}}
    \hspace*{20pt}
    \subfloat[\centering  Time consumption~(orange, right y-axis), Memory usage~(purple, left y-axis)]{\includegraphics[width=0.45\textwidth]{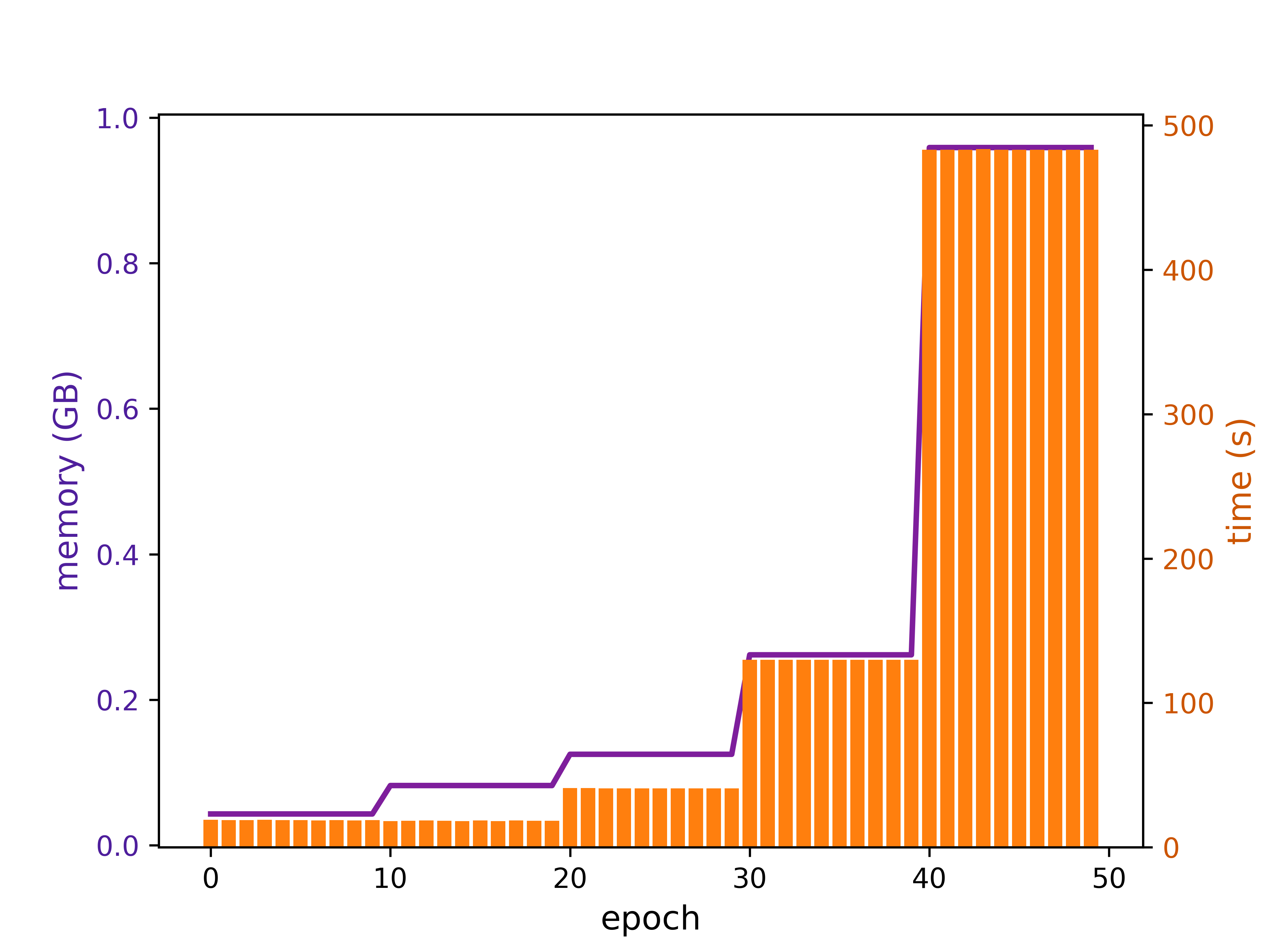}}
    \caption{Training history. The left shows the loss history over batches, and the right shows the time consumption and memory usage over epochs. The loss increases slightly, while the time and memory increase exponentially as the grid size grows after every 10,000 batches (10 epochs).}
    \label{fig:training_history}
\end{figure}

In Fig.~\ref{fig:training_history}, we present the training history of the proposed solver. As per the training settings, the grid size used for training is doubled every 10,000 batches (10 epochs). This change is marked by a minor increase in loss and a significant surge in time consumption and memory usage. The slight uptick in loss, combined with the substantial increase in time and memory, underscores the importance of our training strategy, which involves training the solver on small grid sizes not larger than $511 \times 511$.

\FloatBarrier
\section{Coefficient Distributions}
\label{apd:coeff_dist}

\subsection{GMG Performance over Different Coefficient Distributions}
\label{sec:GMG_results}
As comparisons, we also test the GMG solver over different coefficient distributions and list the number of iterations in Table~\ref{tab:GMG_results}. 
\begin{table}[htbp]
    \centering
    \caption{The number of iterations of the GMG solver over different coefficient distributions.}
    \label{tab:GMG_results}
    \begin{tabular}{|c|l|l|l|l|l|l|l|}
    \hline
        grid & noise & CIFAR10 & FMNIST & MNIST & mldata & Re$10$ & Re$10^5$ \\ \hline
        31x31 & 16.7  & 13 & 16.3  & 19.1  & 13 & 17.6  & 15 \\ 
        63x63 & 19.2 & 15.1 & 16.8 & 18.1 & 15 & 20 & 17.9 \\ 
        127x127 & 21 & 19.8 & 19.5 & 19 & 20 & 21 & 20 \\ 
        255x255 & 25.9 & 25.9 & 24.1 & 23.6 & 26.1 & 24.7 & 26 \\ 
        511x511 & 33.7 & 34.3 & 31.9 & 29.4 & 35.4 & 31 & 34.1 \\ 
        1023x1023 & 45 & 46.3 & 41.8 & 37.4 & 48 & 40 & 46 \\ 
        2047x2047 & 61 & 61.3 & 51.1 & 48 & 65.4 & 53 & 63 \\ 
        4095x4095 & 84 & 85.8 & 71.7 & 61.7 & 87.2 & 71 & 87 \\ \hline
    \end{tabular}
\end{table}

\FloatBarrier
\subsection{Different Coefficient Distributions and Convergence History of the Proposed Solver}
\label{sec:our_solver_results}
\begin{figure}[htbp]
    \centering
    \subfloat[CIFAR10]{\includegraphics[width=0.41\textwidth]{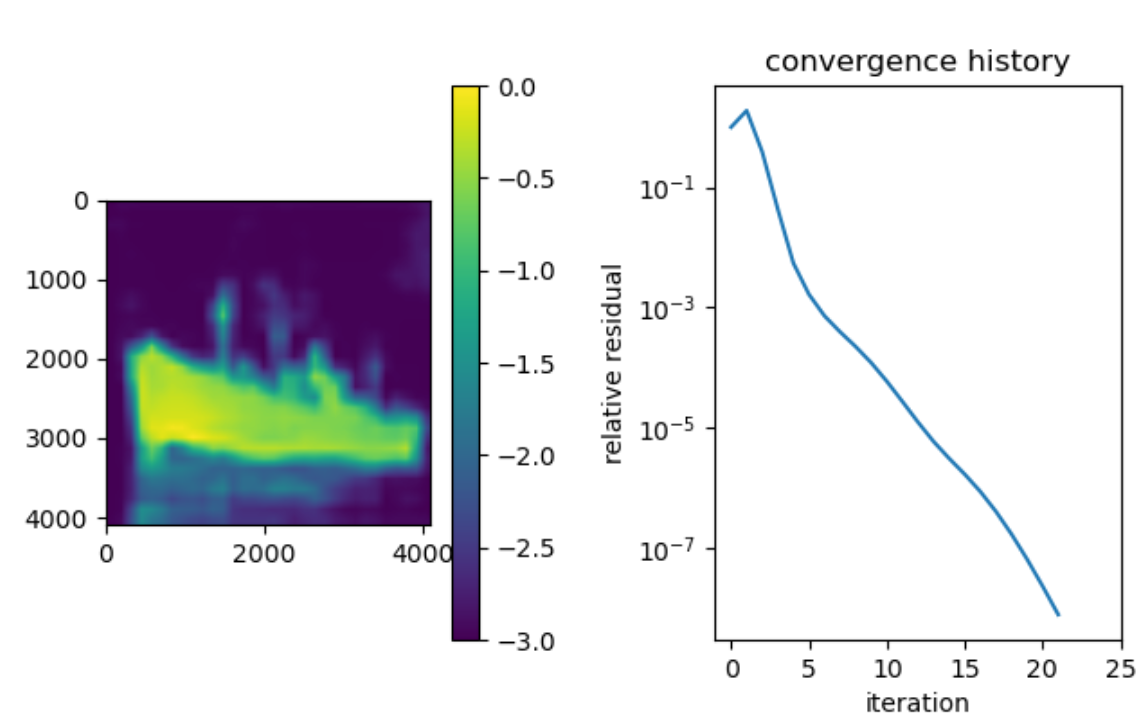}}
    \hspace*{20pt}
    \subfloat[FMNIST]{\includegraphics[width=0.41\textwidth]{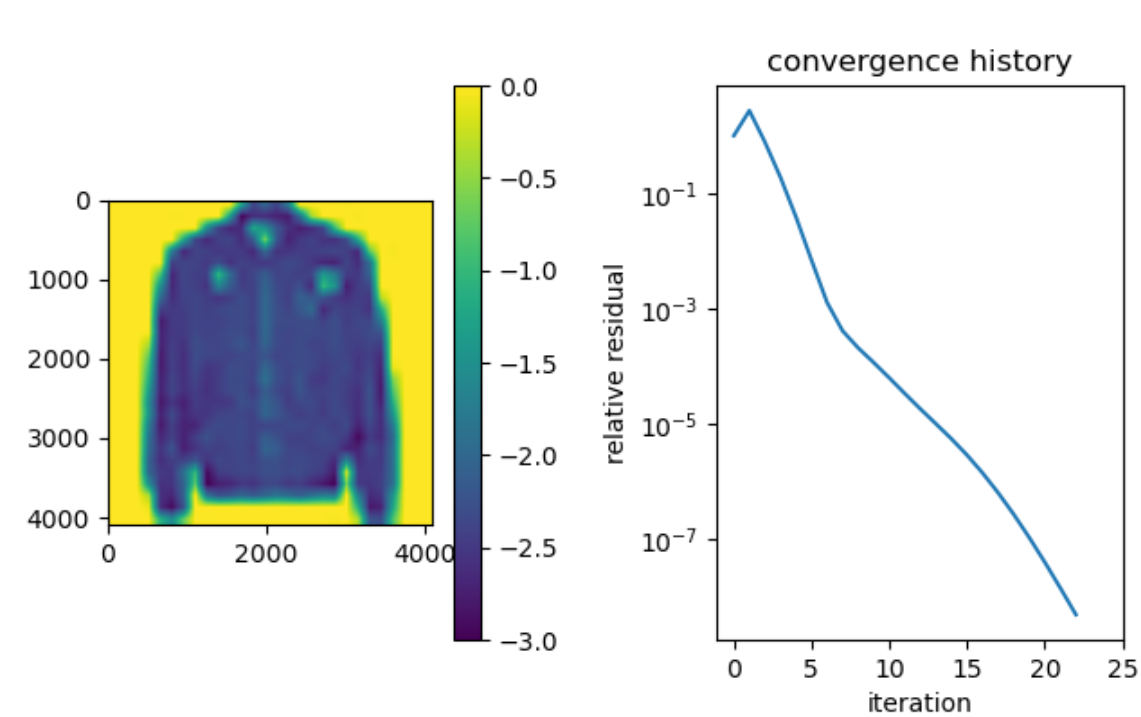}}
    \\
    \subfloat[MNIST]{\includegraphics[width=0.41\textwidth]{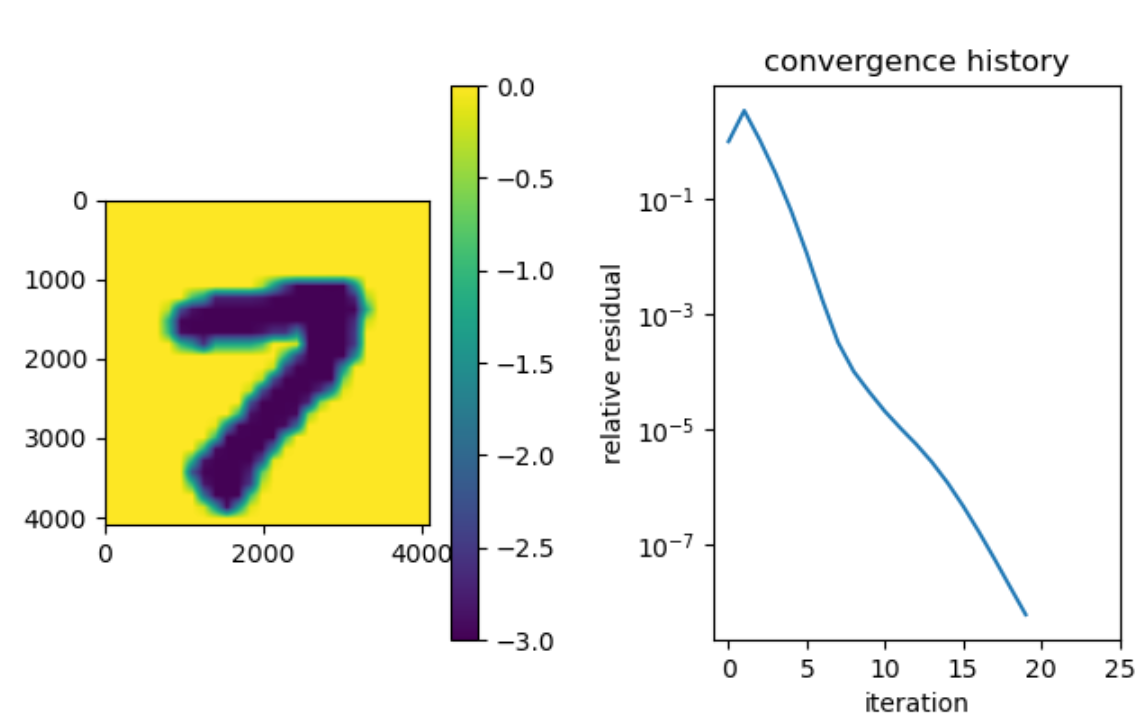}}
    \hspace*{20pt}
    \subfloat[mldata]{\includegraphics[width=0.41\textwidth]{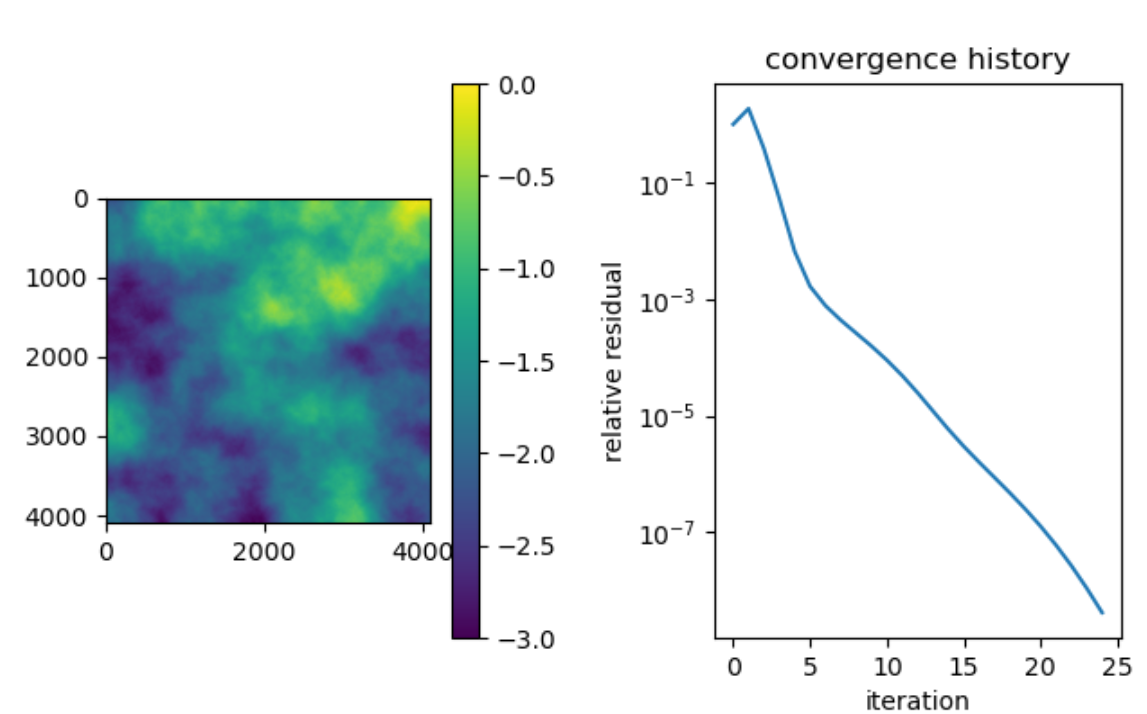}\label{fig:mldata}}
    \caption{The $\log_{10}$ scale coefficient tensor and convergence history. \footnotesize The proposed solver is trained over the mixed coefficient dataset, a combination of white noise and CIFAR10. The right-hand-side is constant, \textbf{1}, and the grid size is $4095 \times 4095$, both of which are not used at training.}
    \label{fig:diff_dist}
\end{figure}

To give a glimpse of different coefficient distributions, we show the $\log_{10}$ scale of coefficient tensors and convergence history of our mixed-CIFAR10-trained solver in Fig.~\ref{fig:diff_dist}. The right-hand-side is constant, \textbf{1}, and the grid size is $4095 \times 4095$, both of which are not used at training. 

\FloatBarrier
\subsection{Algorithm for Generating Multi-Level Data}
\label{sec:mldata}
Despite the white noise distribution and authentic world images, we design a multi-level data generation algorithm to test the impact of data distribution. Basically, we start from a small white noise tensor, interpolate it to a larger size, and add some white noise to it. We repeat this process for several levels and finally interpolate it to the desired size. The detailed algorithm is shown in Alg.~\ref{alg:mldata}. In our experiments, we set $init\_size$ as 5. One can view Fig.~\ref{fig:mldata} to have a glimpse of this generated distribution.
\begin{algorithm}[htbp]
    \caption{Generate Multi-Level Data}
    \label{alg:mldata}
    \begin{algorithmic}[1]
    \Procedure{gen\_data\_multi\_level}{$goal\_size, levels, init\_size$}
        \State $data \gets white\_noise(init\_size)$
        \For{$i$ in $1, 2, ..., levels$}
            \State $data \gets \text{interpolate to } size = 2\times data.size$
            \State $noise\_ratio \gets 2^{-i}$
            \State $data \gets data + noise\_ratio \times white\_noise(data.size)$
        \EndFor
        \State $data \gets \text{interpolate to } size = goal\_size$
    \EndProcedure
    \end{algorithmic}
\end{algorithm}

\FloatBarrier
\bibliographystyle{elsarticle-num}

\bibliography{references}  

\end{document}